\documentclass[a4paper, 10pt]{article}
\usepackage[normalem]{ulem}
\usepackage{cancel}
\usepackage[colorlinks=true]{hyperref}
\hypersetup{urlcolor=blue,citecolor=red}
\usepackage[active]{srcltx}
\usepackage[dvips]{graphicx}
\usepackage{amssymb,amsmath,amsfonts,bbm,pifont,upgreek,bbold,accents} 
\newcommand\beq[1]{ \begin{equation}\label{#1} }
\newcommand{\eeq}{ \end{equation} }
\newcommand{\beqno}{ \[ }
\newcommand{\eeqno}{ \] }
\newcommand\beqa[1]{ \begin{eqnarray} \label{#1}}
\newcommand{\eeqa}{ \end{eqnarray} }
\newcommand{\beqano}{ \begin{eqnarray*} }
\newcommand{\eeqano}{ \end{eqnarray*} }
\newcommand\arr[1]{\left\{\begin{array}{l}#1\end{array}\right.}

\newtheorem{theorem}{Theorem}[section]
\newtheorem{definition}{Definition}[section]
\newtheorem{proposition}{Proposition}[section]
\newtheorem{lemma}{Lemma}[section]
\newtheorem{sublemma}{Sublemma}[section]
\newtheorem{remark}{Remark}[section]
\newtheorem{notationalremark}{Notational Remark}[section]
\newtheorem{corollary}{Corollary}[section]
\newtheorem{assumption}{Assumption}[section]
\newtheorem{claim}{Claim}[section]

\newtheorem{tools}{$\negsp\negsp$}[subsection]

\newcommand\thm[1]{ \begin{theorem}\label{#1}}
\newcommand\thmtwo[2]{ \begin{theorem}[#1]\label{#2}}
\newcommand\ethm{ \end{theorem} }
\newcommand\dfn[1]{ \begin{definition}\label{#1} \rm}
\newcommand\dfntwo[2]{ \begin{definition}[#1]\label{#2} \rm}
\newcommand\edfn{ \end{definition} }
\newcommand\pro[1]{ \begin{proposition}\label{#1}}
\newcommand\protwo[2]{ \begin{proposition}[#1]\label{#2}}
\newcommand\epro{ \end{proposition} }
\newcommand\lem[1]{ \begin{lemma}\label{#1}}
\newcommand\lemtwo[2]{ \begin{lemma}[#1]\label{#2}}
\newcommand\elem{ \end{lemma} }
\newcommand\sublem[1]{ \begin{sublemma}\label{#1}}
\newcommand\sublemtwo[2]{ \begin{sublemma}[#1]\label{#2}}
\newcommand\esublem{ \end{sublemma} }
\newcommand\rem[1]{ \begin{remark}\label{#1} \rm}
\newcommand\erem{ \end{remark} }
\newcommand\notrem[1]{ \begin{notationalremark}\label{#1} \rm}
\newcommand\enotrem{ \end{notationalremark} }
\newcommand\cor[1]{ \begin{corollary}\label{#1}}
\newcommand\cortwo[2]{ \begin{corollary}[#1]\label{#2}}
\newcommand\ecor{ \end{corollary} }
\newcommand\asmp[1]{ \begin{assumption}\label{#1}}
\newcommand\asmptwo[2]{ \begin{assumption}[#1]\label{#2}}
\newcommand\easmp{ \end{assumption} }
\newcommand\clm[1]{ \begin{claim}\label{#1}}
\newcommand\eclm{ \end{claim} }
\newcommand{\proof}{\par\medskip\noindent{\bf Proof}}
\newcommand\equ[1]{{\rm (\ref{#1})}}

\newcommand\ovl[1]{ \overline {#1} }
\newcommand\su[1]{ \frac{1}{ {#1}} }

\newcommand{\torus}{ {\mathbb T}   }

\newcommand{{\real}}{ {\mathbf R}   }
\newcommand{{\integer}}{ {\mathbf Z}   }

\renewcommand{\a }{ {\alpha}   }
\renewcommand{\b}{ {\beta}   }
\newcommand{\g}{ {\gamma}   }

\renewcommand{\d}{ {\delta}   }

\newcommand{\e }{ {\epsilon}   }

\renewcommand{\k}{ {\kappa}   }
\renewcommand{\l}{ {\lambda}   }
\renewcommand{\L}{ {\Lambda}   }
\newcommand{\m}{ {\mu}   }
\newcommand{\n}{ {\nu}   }

\newcommand{\p}{ {\pi}   }
\renewcommand{\P}{ {\Pi}   }
\renewcommand{\r}{ {\rho}   }
\newcommand{\s}{ {\sigma}   }

\renewcommand{\t}{ {\tau}   }
\newcommand{\f}{ {\varphi}   }

\renewcommand{\o}{ {\omega}   }
\renewcommand{\O}{ {\Omega}   }

\newcommand{\const}{{\, \rm const\, }}

\renewcommand{\Im}{{\, \rm Im\, }}
\renewcommand{\Re}{{\, \rm Re\, }}

\newcommand{\cA}{ {\cal A} }
\newcommand{\cB}{ {\cal B} }

\newcommand{\cT}{ {\cal T} }

\newcommand{\cK}{ {\cal K} }
\newcommand{\cC}{ {\cal C} }
\newcommand{\cD}{ {\cal D} }

\newcommand{\cG}{ {\cal G} }
\newcommand{\cL}{ {\cal L} }
\newcommand{\cM}{ {\cal M} }

\newcommand{\cP}{ {\cal P} }

\newcommand{{\cJ}}{ {\cal J} }

\newcommand\ppu{{ (1) }}
\newcommand\ppd{{ (2) }}
\newcommand\ppt{{ (3) }}

\newcommand\ppj{{ (j) }}

\newcommand\ppi{{ {\rm (i)} }}

\newcommand\mm{{\rm m}}
\newcommand\MM{{\rm M}}
\usepackage{color}
\definecolor{applegreen}{rgb}{0.55, 0.71, 0.0}
\definecolor{amber(sae/ece)}{rgb}{1.0, 0.49, 0.0}
\definecolor{amethyst}{rgb}{0.6, 0.4, 0.8}

\newcommand\GG{{\rm G}}

\newcommand\ZZ{{\rm Z}}

\newcommand\PP{{\rm P}}
\newcommand\ii{{\rm i}}
\newcommand\meas{{\, \rm meas}}

\newcommand\id{{\, \rm id \,}}

\newcommand\tk{{ k }}
\newcommand\fg{{\g_1}}
\newcommand\sg{{\ovl\g_2}}
\newcommand\ssg{{\g_2}}
\newcommand\nf{{\rm h}}
\newcommand\nfj{{\rm h}_}
\newcommand\pert{{f}}
\newcommand\so{{\o_2}}
\newcommand\Hs{{M}}
\newcommand\sHs{{\widehat M}}
\newcommand\inHs{{\ovl M}}

\newcommand\HH{{\rm H}}

\newcommand\hh{{\rm h}}
\newcommand{\CC}{{\rm C}}

\begin{document}

\title{Quantitative {\sc kam} theory, with  an application to the three--body problem} 


\author{Gabriella Pinzari\footnote{Department of Mathematics, University of Padova, Padova, Italy. Corresponding Author. {\tt pinzari@math.unipd.it}}, Xiang Liu\footnote{Department of Mathematics, University of Padova, Padova, Italy.}
}


\maketitle 

\begin{abstract}
Based on quantitative ``{\sc kam} theory'', we state and prove two  theorems about the continuation of  maximal and whiskered quasi--periodic motions to slightly perturbed systems exhibiting proper degeneracy. Next, we apply such results to prove that, in the three--body problem, there is a small set in phase space where it is possible to detect both such families of tori. We also estimate the density of such motions in proper ambient spaces.\\
 Up to our knowledge, this is the first proof of co--existence of stable and whiskered tori in a physical system.
\end{abstract}
\tableofcontents

\newpage
\section{Overview}\label{Existence of tori}

\paragraph{1.1 Two {\sc kam} theorems  for properly--degenerate Hamiltonian systems}  We deal with Hamiltonians which meet the demand of being close--to--be--integrable (see, e.g.,~\cite{gallavotti86}), but, in addition, with the number of degrees of freedom of perturbing term  being possibly larger than the one of the unperturbed part. Such kind of Hamiltonians often arise in problems of celestial mechanics, and are 
 referred to as ``properly--degenerate'', after~\cite{arnold63}. We denote them as
$$H(I, \varphi,  p, q;   \m)=H_0(I)+\m\,P(I, \f, p, q;   \m)\,,$$
where the coordinates $(I, \varphi)=(I_1, \ldots, I_n, \varphi_1, \ldots, \varphi_n)$ are of ``action--angle'' kind (after a possible application of the Liouville--Arnold theorem to the unperturbed term), while (for our needs) the $(p, q)=(p_1, \ldots, p_m, q_1, \ldots, q_m)$ are ``rectangular'', namely,
take value in a small ball (say, of radius $\varepsilon_0$) about some point (say, the origin). The
symplectic form is standard:
$$\Omega=dI\wedge d\varphi+dp\wedge dq=\sum_{i=1}^n dI_i\wedge d\varphi_i+\sum_{i=1}^m dp_i\wedge dq_i\,.$$ 
We work in the real--analytic framework, which means that we assume that $H$ admits a holomorphic extension on a complex neighburhood of the real ``phase space'' (namely, the domain) $$\cP_{\varepsilon_0}:=V\times {\mathbb T}^n\times B^{2m}_{\varepsilon_0}\,,$$ where $V\subset {\mathbb R}^n$ is bounded, open and  connected, ${\mathbb T}={\mathbb R}/(2\pi {\mathbb Z})$ is the ``flat torus'', $B^{2m}_{\varepsilon_0}$ is the $2m$--dimensional ball around $0$ of radius ${\varepsilon_0}$, relatively to some norm in $\mathbb R^{2m}$.

\noindent
In this framework, we present\footnote{We refer to specialized literature for hystorical notices and constructive approaches to {\sc kam} theory: see, e.g.,~\cite{gallavotti1994, gallavottiG1995, bonettoGGM1998, chierchiaPr2019} and references therein.} two ``{\sc kam} theorems''  which deal with different situations. 
A basic assumption, common to both statements, and  often referred to as ``Kolmogorov condition'', is:
\begin{itemize}
\item[\rm (A${}_1$)]  
the map $I\to \partial_I H_0(I)$ is a diffeomorphism of $V$.
\end{itemize}
However, due to the proper degeneracy mentioned above, such assumption is to be reinforced with some statement concerning the perturbing term, or, more precisely, its Lagrange average 
$$P_{\rm av}(I, p, q; \mu):=\frac{1}{(2\pi)^n}\int_{[0, 2\pi]^n}P(I, \varphi, p, q; \mu)d^n\varphi$$
with respect to the $\varphi$--coordinates. Such extra--assumption will be different in the two statements, therefore we quote them below.

\vskip.1in
\noindent
 The first result is  a revisitation of the so--called ``Fundamental Theorem'' by V.I. Arnold,~\cite{arnold63}. Such theorem has been already studied, generalised and extended in previous works~\cite{chierchiaPi10, pinzari18a}. Here we deal 
 with the situation (not considered in the aforementioned papers) where $P_{\rm av}$ admits a ``Birkhoff Normal Form'' ({\sc bnf} hereafter) about $(p, q)=(0, 0)$ of high\footnote{\cite{arnold63, chierchiaPi10, pinzari18a} deal with the ``minimal'' case $s=2$. The case $s=2$ is called here ``minimal'' as we work in the framework of generalizations of the Kolmogorov condition $(A_1)$ above. In~\cite{Russmann:2001, fejoz04}, using different techniques, the case $s=1$ has been considered.} order; say $s$.  As expected, a higher order of {\sc bnf}  allows to  improve the measure of the ``Kolmogorov set'', namely the set given by the union of all {\sc kam} tori. We shall prove\footnote{For simplicity of notations, we do not write  $\mu$ among the arguments of the functions in Theorem~\ref{stable toriREF} and~\ref{thm:simplifiedHYPER}.} the following

\begin{theorem}\label{stable toriREF}
Assume $(A_1)$ above and the following conditions:
\begin{itemize}

\item[\rm (A${}_2$)]  $ P_{\rm av} (I, p, q)=\sum_{j=1}^{s}\cP_j(r;I)
+{\rm O}_{2 s+1}(p,q;I)$,  with  $r_i:=\frac{p_i^2+q_i^2}{2}$ and $\cP_j(r;I)$ being a polynomial of degree $j$ in $r=(r_1, \cdots, r_{m})$, for some $2\le s\in {\mathbb N}$.

\item[\rm (A${}_3$)] the $m\times m$ matrix $\beta (I)$ of the coefficients of the second--order term $\cP_2(r;I)=\frac{1}{2}\sum_{i,j=1}^{m} \b_{ij}(I)r_i r_j$ is non degenerate:
 $|\det \b(I)|\ge \const>0$ for all $I\in V$.
\end{itemize}
Then, there exist positive numbers $\varepsilon_*<\varepsilon_0$, $C_*$ and $c_*$ such that, for 
\beqa{epsmuREF}0<\varepsilon<\varepsilon_*\ , \qquad 
0<\m<\frac{\varepsilon^{2s+2}}{C_* (\log \varepsilon^{-1})^{c_*}}\ .
\eeqa
one can find a 
set $\cK\subset  \cP_\varepsilon$  formed by the union of $H$--invariant $n$--dimensional Lagrangian  tori, on which the $H$--motion is analytically conjugated to  linear Diophantine quasi--periodic motions with frequencies $(\o_1,\o_2)\in{\mathbb R}^{n_1}\times{\mathbb R}^{n_2}$ with $\o_1=O(1)$ and $\o_2=O({\m})$. The  set $\cK$  has positive Liouville--Lebesgue measure and satisfies
\beqa{meas est*REF}
\meas \cP_\varepsilon>\meas \cK>  \Big(1- C_* {\varepsilon}^{s-\frac{3}{2}}\Big)\meas \cP_\varepsilon\ .
\eeqa
\end{theorem}

\vskip.1in
\noindent
The second result deals with lower--dimensional quasi--periodic motions, the so--called ``whiskered tori''. These are  $n$--dimensional
quasi--periodic motions (in a phase space of dimension $2n+2m$),  approached or reached at an exponential rate.
For simplicity, and in view of our application, we focus on the case $m=1$. In addition, we allow a further degeneracy in the Hamiltonian: the unperturbed term $H_0$ may possibly depend not on all the $I$'s, but only on a part of them.

\begin{theorem}\label{thm:simplifiedHYPER}
Let $m=1$, and let $H_0$ depend on the components $I_1=(I_{11}, \ldots, I_{1n_1})$  of the $I=(I_1, I_2)$'s, with $1\le n_1\le n:=n_1+n_2$. Assume $(A_1)$ with $I_1$ replacing $I$ and, in addition,
 that 
\begin{itemize}
\item[\rm (A${}_2'$)]  
$ P_{\rm av} (I, p, q;   \m)=P_0 (I, pq;   \m)+ P_1 (I,\varphi, p, q;   \m)$ 
  with $\|P_1\|\le  a \| P_0\|$;

\item[\rm (A${}_3'$)]  
 $| \partial_{pq} P_0 |\ge \const>0$ {and $|\det \partial^2_{I_2} P_0 |\ge \const>0$ if $n_2\ne 0$}.
\end{itemize}
Fix $\eta>0$.
Then, there exist positive numbers $a_*$, $\varepsilon_*<\varepsilon_0$, $C_*$ and $c_*$ such that, if 
\beqa{epsmuHYPER}
0<\varepsilon<\varepsilon_*\,,\quad 0< a < a _*\varepsilon^4\ ,\quad 
0<\m<\frac{C_* (a \|P_0\|)^{1+\eta}}{ (\log  a ^{-1})^{c_*}}
\eeqa
 one can find a 
set $\cK$  formed by the union of $H$--invariant $n$--dimensional Lagrangian  tori, on which the $H$--motion is analytically conjugated to  linear Diophantine quasi--periodic motions with frequencies $(\o_1,\o_2)\in{\mathbb R}^{n_1}\times{\mathbb R}^{n_2}$ with $\o_1=O(1)$ and $\o_2=O({\m})$. The projection $\cK_0$ of set $\cK$  on $\cP_0:=V\times {\mathbb T}^{n}$ has positive Liouville--Lebesgue measure and satisfies
\beqa{meas est*HYPER}
\meas {\cal P}_0>\meas \cK_0> \Big(1- C_* \sqrt{ a }\Big) \meas {\cal P}_0\ .
\eeqa
Furthermore, for any $\cT\in \cK$ there exist two $(n+1)$--dimensional invariant manifolds ${\cal W}_{\rm u}$, ${\cal W}_{\rm s}\subset  \cP_{\varepsilon_*}$  such that $\cT={\cal W}_{\rm u}\cap {\cal W}_{\rm s}$ and the motions on ${\cal W}_{\rm u}$, ${\cal W}_{\rm s}$ leave, approach $\cT$ at an exponential rate.
\end{theorem}
Before we go on with describing how we aim to use the theorems above, we premise some comment.
\begin{itemize}\rm 
\item[(i)] The conditions involving $\mu$ in~\equ{epsmuREF} and~\equ{epsmuHYPER} are not optimal. With a procedure similar to the one shown in~\cite[proof of Theorem 1.2, steps 1--4]{chierchiaPi10}, one can show that they can be  relaxed to, respectively
$$\m<\su{C_* (\log \varepsilon^{-1})^{2b}}\,,\qquad \m<\su{C_* (\log (a \|P_0\|)^{-1})^{2b}}$$
with some  $C_*$, $b>0$.
\item[(ii)] The careful bounds on the measure of the invariant sets provided in~\equ{meas est*REF}  and~\equ{meas est*HYPER}  are needed in view of our application. 
 Indeed, we shall apply both the theorems above in order to prove that, in the three--body problem, closely to the  {\it co--planar, co--circular, outer retrograde configuration} (see below for the exact definition),  full--dimensional and ``whiskered'' quasi--periodic tori co--exist (the result was conjectured in~\cite{pinzari18}). In the application, $\varepsilon$ will correspond  to the maximum eccentricity or inclination; $a$ the semi--major axes ratio, and the use of a high--order {\sc bnf} in Theorem~\ref{stable toriREF} will be necessary because the size of the set goes to $0$ with some power of $\varepsilon$ ($s=4$ will be enough for our application).  
 \item[(iii)] Following~\cite{chierchiaG94}, Theorem~\ref{thm:simplifiedHYPER} might be extended to prove the existence of ``diffusion paths'' and ``whisker ladders''. We shall not do, as proving Arnold instability (in the sense of~\cite{arnold64}) for the system~\equ{helio} below is not the purpose of this paper. We however remark that such kind of instability has been found for  the {\it four--body problem} in a very similar framework~\cite{clarkeFG22}.
We remark that proofs of chaos or Arnold instability in celestial mechanics are quite recent 
~\cite{FGKR14, delshamsKDRS2019}, {by the difficulty of overcoming the so--called ``problem of large gaps''. See \cite{guzzoEP2020} and references therein.}
\item[(iv)] Another important aspect in view of the application described above
 is a rather standard consequence of the proof of Theorem~\ref{thm:simplifiedHYPER}: If  $P$ (namely, $P_1$) has an equilibrium at $(p, q)=0$, then, along the motions of $\cal K$, the coordinates $(p, q)$ remain fixed at $(0, 0)$ (rather than varying closely to it), namely
 $${\cal K}\subset V\times{\mathbb T}^{n}\times\{(0, 0)\}\,.$$
 More generally, the stable and unstable invariant manifolds do not shift from the unperturbed ones:
$${\cal W}_{\rm s}\subset\cP_{\varepsilon}\cap \big\{q=0\big\}\,,\qquad {\cal W}_{\rm u}\subset \cP_{\varepsilon}\cap \big\{p=0\big\}\,.$$
\end{itemize}

\paragraph{1.2 Application to the three--body problem} 
We apply the results above to prove that, in a region of the phase space of the three--body problem, and under conditions that will be specified later, full dimensional and whiskered tori co--exist. We underline that the co--existence of such different kind of motions is not a mere consequence of the non--integrability of the system (as in such case the result would be somewhat expected) as it persists in two suitable 
integrable approximations of the system, close one to the other. Indeed, such motions will be found in a very small zone in the phase space of the three--body problem which simultaneously is in the neighborhood of an elliptic equilibrium of one of such approximations and in a hyperbolic one of the other. Such an occurrence is intimately related to the use of two different systems of coordinates, which are singular one with respect to the other, in the region of interest. The authors are not aware of the appearance of such phenomenon, previously.

\vskip.1in
\noindent
After the ``heliocentric reduction'' of translational invariance, the three--body problem Hamiltonian with gravitational masses equal to $m_0$, $\m m_1$ and $\m m_2$ and Newton constant $\cG\equiv 1$, takes the form of the two--particle system (see, e.g,~\cite{fejoz04, laskarR95} for a derivation):
\beqa{helio}\HH_{\rm 3b}=\sum_{i=1}^2\left(\frac{|y^\ppi|^2}{2{\rm  m}_i}-\frac{{\rm  m}_i{\rm  M}_i}{|x^\ppi|}\right) +\mu\left(-\frac{m_1m_2}{|x^\ppu-x^\ppd|}+\frac{y^\ppu\cdot y^\ppd}{m_0}\right)
\eeqa
with suitable values of $\mm_i=m_i+{\rm O}(\mu)$, $\MM_i=m_0+{\rm O}(\mu)$. We consider the system in the Euclidean space, namely we take, in~\equ{helio}, $y^\ppi$, $x^\ppi\in \mathbb R^3$, with $x^\ppu\ne x^\ppd$.\\
We call {\it Kepler maps}
the class of symplectic\footnote{Namely
 verifying 
$$\Omega=d\L\wedge d\ell+d{\rm y}\wedge d{\rm x}=\sum_{i=1}^2 d\L_i\wedge d\ell_i+\sum_{i=1}^4d{\rm y}_i\wedge d{\rm x}_i\,.$$}
 coordinate systems $\cC=(\L_1, \L_2, \ell_1, \ell_2,  \rm y, \rm x)$ for the Hamiltonian~\equ{helio}, where
$ \rm y=(y_1, \ldots, y_4), \rm x=(x_1, \ldots, x_{4})$,
such that:

\begin{itemize}
\item[--] $\Lambda_i=\mm_i \sqrt{\MM_i a_i}$, where $a_i$ denotes the  semi--major axis of the $i^{\rm th}$ instantaneous\footnote{With reference to the three--body  Hamiltonian~\equ{helio},
 the
{\it $i^{\rm th}$ instantaneous ellipse} is the orbit generated by
$\hh_i:=\frac{|y^\ppi|^2}{2{\rm  m}_i}-\frac{{\rm  m}_i{\rm  M}_i}{|x^\ppi|}$
in a region of phase space where $\hh_i$ is negative.} ellipse; 
\item[--]  $\ell_1, \ell_2\in \mathbb T$ are  conjugated to $\L_1$, $\L_2$. Such angles are defined in a different way according to the choice of $\cC$. In all known examples, they are related to the area spanned by the planet along the instantaneous ellipse.
\end{itemize}  
Using a Kepler map, 
the Hamiltonian~\equ{helio} takes the form
\beqa{HC}\HH_\cC=-\frac{\mm^3_1\MM^2_1}{2\L_1^2}-\frac{\mm^3_2\MM^2_2}{2\L_2^2}+\m f_\cC(\L_1, \L_2, \ell_1, \ell_2,  {\rm \hat y, \rm \hat x})\eeqa
where
$ {\rm \hat y, \rm \hat x}$ include the couples  $ ({\rm y}_i, {\rm x}_i)$ such\footnote{The reason of this is that the Hamiltonian~\equ{helio} has  first integrals, as recalled in the next section.} that nor ${\rm y}_i$ nor ${\rm x}_i$ is negligible. $ {\rm \hat y, \rm \hat x}$ are often called {\it degenerate coordinates}, because  they do not appear in~\equ{HC} when $\mu$ is set to zero.
In other words, $\HH_\cC$ is a properly--degenerate close--to--be--integrable system, in the sense of the previous paragraph.
\vskip.05in
\noindent
We call {\it co--planar, co--circular, outer retrograde configuration} the configuration of two planets in circular and 
co--planar motions, with the  the angular momentum of the outer planet having opposite verse to the resulting one.
In~\cite{pinzari18} it has been pointed out that, under a careful choice of $\cC$ such configuration plays the r\^ole of an equilibrium 
for the $(\ell_1, \ell_2)$--averaged perturbing function
$$\ovl f_\cC(\L_1, \L_2,  {\rm \hat y, \rm \hat x})=\frac{1}{(2\pi)^2}\int_{[0, 2\pi]^2}f_\cC(\L_1, \L_2, \ell_1, \ell_2,  {\rm \hat y, \rm \hat x}) d\ell_1d\ell_2\,.$$
 But what matters more is that, closely to such equilibrium, there exist two such  $\cC_i$'s such that
the Hamiltonian $\HH_{\cC_1}$ is suited to Theorem~\ref{stable toriREF}, while $\HH_{\cC_2}$ is suited to Theorem~\ref{thm:simplifiedHYPER}. This leads to the following result, which states co-existence of stable and whiskered quasi--periodic motions in the three--body problem. It will be made more precise (see Theorem~\ref{main}  below)  and proved along the paper.

\vskip.1in
\noindent
{\bf Theorem A} {\it
In {the vicinity} of the co--planar, co--circular, outer retrograde configuration, and provided that the masses of the planets and the semi--axes ratio are small, there exists a positive measure set $\cK_1$ made of $5$--dimensional quasi--periodic motions $\cT_1$'s ``surrounding'' $($in a sense which will be specified$)$ $3$--dimensional 
quasi--periodic motions $\cT_2$'s, each  equipped with two invariant manifolds, called, respectively, unstable, stable manifold, where the motions  are respectively asymptotic  to the $\cT_2$'s in the past, in the future.
}
\vskip.1in
\noindent
We conclude with saying how this paper is organized. 
\begin{itemize}
\item[\textbullet] In Sections~\ref{The elliptic character}  and~\ref{The hyperbolic character} we recall the main arguments of the discussion in~\cite{pinzari18}, which lead to put the system~\equ{helio} to a form suited to apply Theorems ~\ref{stable toriREF} and~\ref{thm:simplifiedHYPER}. 
\item[\textbullet] In Sections~\ref{Existence and co--existence of two families of tori} and~\ref{Proof of Propositions} we check that the two domains where Theorems ~\ref{stable toriREF} and~\ref{thm:simplifiedHYPER} apply have a non--empty intersection, and such intersection includes both families of tori. This check is subtle, because of the difference of the frameworks used. 
\item[\textbullet] In Section~\ref{KAM theory} we prove Theorems ~\ref{stable toriREF} and~\ref{thm:simplifiedHYPER} via a carefully quantified {\sc kam} theory.
\end{itemize}

\section{Ellipticity and hyperbolicity closely to co--planar, co--circular, outer retrograde configuration}\label{Double nature of co--planar, co--circular, outer retrograde configuration}

Putting the system in a form suited to Theorem ~\ref{stable toriREF} requires identifying an elliptic equilibrium; while Theorem~\ref{thm:simplifiedHYPER} calls for a hyperbolic one. \\
Denoting as
$\CC^{(j)}:=x^{(j)}\times y^{(j)}$ the angular momenta of the planets, we
proceed to study motions evolving from initial data close to the manifold
    \beqa{singularities1}
{\cal M}_\p:=\Big\{(y, x):\ 
\CC^\ppu\parallel(-\CC^\ppd)\parallel\CC
\,,\ {\rm and}\
x^\ppu\,,\ x^\ppd \ {\rm describe\ circular\ motions
\,.}
\Big\}\,.
\eeqa 
The sub--fix ``${\pi}$'' recalls that $\rm C^\ppu$
and $\rm C^\ppd$ are  opposite. 
In the two next sections, we recall material from~\cite{pinzari18}, which highlights
a sort of ``double (elliptic, hyperbolic) nature'' of ${\cal M}_\p$.

\subsection{Ellipticity (with  {\sc bnf})}\label{The elliptic character}  
Basically\footnote{As pointed out in~\cite{pinzari18}, the only note--worthing difference with the case studied  in~\cite{chierchiaPi11b} (which deals with prograde motions of the  planets, namely, revolving all in the same verse) is that here the elliptic  character of the equilbrium  does not follow for free from the symmetry of the Hamiltonian, but is checked manually. }, the construction of the elliptic equilibrium -- and of its associated {\sc bnf} -- proceeds as in~\cite{chierchiaPi11b}. We briefly resume the procedure here.

\vskip.1in
\noindent
We fix a domain $\cD_{{\textrm{\sc c}}}\subset {\mathbb R}^{12}$ for {impulse--position ``Cartesian'' coordinates} $$\textrm{\sc c}=(y, x):=(y^\ppu, y^\ppd, x^\ppu, x^\ppd)$$ of two point masses  relatively to a prefixed orthonormal frame $(k^\ppu, k^\ppd, k^\ppt)$ in ${\mathbb R}^3$. 
As a first step, we switch to a set of coordinates, well known in the literature, which we name {\sc jrd}, after C. G. J. Jacobi, R. Radau and A. Deprit ~\cite{jacobi1842, radau1868, deprit83}, who at different stages, contributed to their construction.

\noindent
We fix a region of phase space where 
the orbits $t\to (x^\ppj(t), y^\ppj(t))$ generated by the  unperturbed ``Kepler'' Hamiltonians
  \beqano
\hh^\ppj_{{\rm  k}}:=\frac{|y^\ppj|^2}{2{\rm  m}_j}-\frac{{\rm  m}_j{\rm  M}_j}{|x^\ppj|} \eeqano
in~\equ{helio}
 are  ellipses with
 non-vanishing eccentricity.  Then we denote as
 ${\rm P}^{(j)}$  the unit vectors pointing in the directions of the perihelia; as $a_j$ the semi--major axes; 
 {as $\ell_j$ the ``mean anomaly'' of $x^\ppj$(which, we recall, is defined as area of the elliptic sector from ${\rm P}^{(j)}$ to $x^\ppj$ ``normalized at $2\p$'');}
 as  ${\rm C}^\ppj=x^\ppj\times y^\ppj$, $j=1$, $2$,  the angular momenta of the two planets and ${\rm C}:={\rm C}^\ppu+{\rm C}^\ppd$ the total angular momentum integral. We assume that the ``nodes''
\beqa{Dep nodes}\n_1:=k^\ppt\times{\rm C}\ ,\quad \n:={\rm C}\times {\rm C}^\ppu={\CC}^\ppd\times{\CC}^\ppu\eeqa
do not vanish, anytime. Such condition is equivalent to ask that
the planes determined by the instantaneous  ellipses and the $(k^\ppu, k^\ppd)$ plane never pairwise coincide. As in previous works, we use the following notations. 
For three vectors $u$, $v$, $w$ with $u$, $v\perp$ $w$, we denote as $\a_{w}(u,v)$  the angle formed by $u$ to $v$ relatively to the positive (counterclockwise) orientation established by $w$. Then the {\sc jrd} coordinates are here denoted with  the symbols
\beq{J}\textrm{\textrm{\sc jrd}}:=\Big(\widehat{\textrm{\sc jrd}}:=(\L_1,\L_2,\GG_1,\GG_2, \ell_1,\ell_2,\g_1,\g_2), (\GG,\ZZ,\g, \zeta)\Big)\in {\mathbb R}^4\times {\mathbb T}^4\times {\mathbb R}^2\times {\mathbb T}^2\eeq
and defined via the formulae
\beqa{coordinates}
\begin{array}{llllrrr}
 \left\{
\begin{array}{lrrr}
 \ZZ:={\rm C}\cdot k^\ppt \\
{\rm G}:=\|{\rm C}\| \\
{\rm G}_1:=\|{\rm C}^\ppu\|\\
{\rm G}_2:=\|{\rm C}^\ppd\|\\
\L_j:={\rm  M}_j\sqrt{{\rm  m}_j a_j}
\end{array}
\right.\qquad
\left\{
\begin{array}{lrrr}
\zeta:=\a_{k^\ppt}(k^\ppu, \n_1)\qquad& \\
 \g:=\a_{{\rm C}}(\n_1, \n)\qquad&\\
 {\g}_1:=\a_{{{\rm C}^\ppu}}(\n, {\rm P}^\ppu)&\\
 {\g}_2:=\a_{{\rm C}^\ppd}(\n,{\rm P}^\ppd)&\\
\ell_j := {\rm mean\ anomaly\ of}\ x^\ppj& 
\end{array}
\right.
\end{array}
\eeqa
The main point of {\sc jrd} is that $\ZZ$, $\zeta$ and $\gamma$ are ignorable coordinates and $\GG$ is constant along the motions of SO(3)--invariant systems. Therefore, most of motions of SO(3)--invariant systems are effectively described by the ``reduced'' coordinates
$\widehat{\textrm{\sc jrd}}$. This strong property cannot be  exploited in the  case study of the paper, as the manifold ${\cal M}_\p$ in~\equ{singularities1} is a singularity of the change~\equ{coordinates}. More generally, any
 co--planar or circular\footnote{Circular configurations correspond to $\GG_i=\Lambda_i$; co--planar configurations correspond to $\GG=\s_1\GG_1+\s_2 \GG_2$, with $(\s_1, \s_2)\in\{\pm 1\}^2\setminus\{(-1, -1)\}$. }
configuration is so.
 Pretty similarly as in~\cite{chierchiaPi11b}, we bypass such difficulty switching to new coordinates denoted as
\beqano
{\textrm{\sc rps}}_{\pi}:=\Big(\widehat{\textrm{\sc rps}}_{\pi}:=(\L_1,\L_2,\l_1,\l_2,   \eta_1,\eta_2,  \xi_1,\xi_2,  p, q),  (Z, \zeta)
\Big) 
\eeqano
where the  $\Lambda_j$'s, $\ZZ$ and $\zeta$ are the same\footnote{There is an inessential difference between the definition~\equ{PR} and the one in~\cite[Eqs. (25), (26))]{pinzari18a}. Denoting as $\ovl{\textrm{\sc rps}}_{\pi}:=\Big((\L_1,\L_2,\ovl\l_1,\ovl\l_2,   \ovl\eta_1,\ovl\eta_2,  \ovl\xi_1,\ovl\xi_2,  \ovl p, \ovl q),  (P,  Q)
\Big)$ the coordinates defined in~\cite{pinzari18a}, we have the following relations:
 {\scriptsize \beqano
\arr{
\displaystyle\ovl\lambda_1=\lambda_1+\zeta
 \\
\displaystyle\ovl\lambda_2=\lambda_2-\zeta
 \\
\displaystyle \ovl\eta_1+\ii\ovl\xi_1=(\eta_1+\ii\xi_1)e^{-\ii\zeta
}\\
\displaystyle \ovl\eta_2+\ii\ovl\xi_2=(\eta_2+\ii\xi_2)e^{\ii\zeta
}\\
\displaystyle \ovl p+\ii \ovl q=(p+\ii q) e^{\ii\zeta
}\\
\displaystyle P+\ii Q=\sqrt{2(\GG-\ZZ)}\,e^{-\ii\zeta}
}
\eeqano}
But as $(Z, \zeta)$, $(P, Q)$ and $\zeta$ do not appear in the Hamiltonian, its  expression does not change.
} as in~\equ{J}, while 
  \beqa{PR}
\arr{
\displaystyle\lambda_1=\ell_1+\gamma_1+\gamma
 \\
\displaystyle\lambda_2=\ell_2+\gamma_2-\gamma
 \\
\displaystyle \eta_1+\ii\xi_1=\sqrt{2(\L_1-\GG_1)}e^{-\ii(\g_1+\g
)}\\
\displaystyle \eta_2+\ii\xi_2=-\sqrt{2(\L_2-\GG_2)}e^{\ii(-\g_2+\g
)}\\
\displaystyle p+\ii q=-\sqrt{2(\GG+\GG_2-\GG_1)} e^{\ii\g
}
}
\eeqa
As in {\sc jrd}, $(Z, \zeta)$ is a  cyclic couple
in SO(3)--invariant Hamiltonians but now no more cyclic coordinates but it appears. This leaves the system with 5 degrees of freedom and an extra--integral: the action $\GG$ written using $\textrm{\sc rps}_\pi$:
\beqa{Grps}{\rm G}_{\textrm{\sc rps}_\pi}:=\Lambda_1-\Lambda_2-\frac{\eta_1^2+\xi_1^2}{2}+\frac{\eta_2^2+\xi_2^2}{2}+\frac{p^2+q^2}{2}\,.\eeqa
 We denote as
\beqa{rps Ham real}
\HH_{{\textrm{\sc rps}}_{\pi}}&:=&-\frac{\mm_1^3\MM_1^2}{2\L_1^2}-\frac{\mm_2^3\MM_2^2}{2\L_2^2}+\m\, \Big(
-\frac{m_1 m_2}{|x_{{\textrm{\sc rps}}_{\pi}}^\ppu-x_{{\textrm{\sc rps}}_{\pi}}^\ppd|}+ \frac{y_{{\textrm{\sc rps}}_{\pi}}^\ppu\cdot y_{{\textrm{\sc rps}}_{\pi}}^\ppd}{m_0}
 \Big)\nonumber\\
 &=:&\hh_{\rm k}(\L)+\m f_{\textrm{\sc rps}_{\pi}}(\L,\l, z)\qquad z:=(\eta,\xi,p,q)\eeqa
the Hamiltonian~\equ{helio} written in  {\sc rps}${}_{\pi}$ coordinates, and worry  about it.

\noindent
We note that the manifold ${\cal M}_\pi$ in~\equ{singularities1} is now given by
$${\cal M}_\pi=\big\{{{\textrm{\sc rps}}_{\pi}}:\ z=0\big\}\,.$$
Then we consider a neighborhood of ${\cal M}_\pi$ 
of the form
 \beqano{\cal M}_{\textrm{\sc rps}_\pi, \varepsilon_0}:={\cal L}\times {\mathbb T}^2\times B^{6}_{\varepsilon_0}(0)\ ,\eeqano
 where  $B^6_{\varepsilon_0}$ is the 6--ball centred at $0\in {\mathbb R}^6$ with radius $\varepsilon_0$; ${\mathbb T}:={\mathbb R}/(2\p{\mathbb Z})$ and ${\cal L}$ is
defined as
  \beq{L0}{\cal L}:=\Big\{\L=(\L_1,\L_2): \ \L_-< \L_2< \L_+\ ,\quad k_-\L_2< \L_1< k_+\L_2\Big\}\,.\eeq
Here, $0<\L_-<\L_+$
 are arbitrarily taken (more conditions on such numbers will be specified in the course of the paper) and,
 for fixed positive\footnote{Observe that $\a_-$ and $\a_+$ have the meaning of lower and upper bound for the semi--major axes ratio $\alpha=a_1/a_2$, namely,
 $$\a_-\le \a\le \a_+\qquad \forall\ (\L_1, \L_2)\in \cL\ .$$
 Indeed, from the formula  \beqano
 \frac{\Lambda_1}{\Lambda_2}=\frac{m_1}{m_2}\sqrt{\frac{m_0+\m m_2}{m_0+\m m_1}\alpha}\eeqano
 we find
 $$\a\le \left(\frac{m_2}{m_1}\right)^2\frac{m_0+\m m_1}{m_0+\m m_2} k_+^2=\a_+$$
 and, similarly, $\a\ge\a_-$.} numbers $0<\a_-<\a_+<1$,
$k_\pm$ are constants depending on $\alpha_\pm$ and the masses via 
 \beq{kpm}k_\pm:=\frac{m_1}{m_2}\sqrt{\frac{m_0+\m m_2}{m_0+\m m_1}\a_\pm}\ .\eeq

\noindent
We now take $0<\delta<1$ and\footnote{The reader might ask the reason of inequalities in~\equ{domain rps+4}.
This is related to the fact that we want to investigate a region of phase space where the inner planet, labeled as ``1'', has a larger angular momentum, namely, $\rm G_1>\rm G_2$, and, simultaneously, the masses of the planets, as well as their eccentricities and mutual inclination are small.  As, when eccentricities and mutual inclination go to zero, the $\rm G_i$ reduce to $\Lambda_i$, by~\equ{L0}, the number $k_-$ in ~\equ{kpm} needs to be strictly larger than $1$. Conditions  ~\equ{domain rps+4} are apt to ensure this, as in fact they immediately imply
$$1-\d \le \frac{m_0+\m m_2}{m_0+\m m_1}\le1+\d$$
hence, by~\equ{kpm},
$$k_-\ge \frac{m_1}{m_2}\sqrt{(1-\d)\a_-}>1\,.$$
} and assume
 \beqa{domain rps+4} 0<\frac{m_2}{m_1}<\min\left\{\sqrt{(1-\d)\a_-}\ ,\ 1-\d\right\}\,,\quad 0<\m<\m_0(\d):=\frac{\d m_0}{m_1(1-\d)-m_2}\eeqa
Then we\footnote{The proof in~\cite[Appendix A]{pinzari18} is given with $\delta= 1-\frac{1}{4\chi^2}\ge\frac{3}{4}$, but works well also for any $\delta\in (0, 1)$. Indeed, for $(\L_1, \L_2)\in {\cal L}$, 
 $$\L_1-\L_2=\L_2\left(\frac{\L_1}{\L_2}-1\right)\ge \L_- (k_--1)\ge \L_- \left(\frac{m_1}{m_2}\sqrt{(1-\d)\a_-}-1\right)\ .$$
 Therefore, for $(\L_1, \L_2)$ on a complex neighborhood of ${\cal L}$ depending on $\L_-$, $m_1$, $m_2$, $\a_-$ and $\d$ we shall have
 $|\L_1-\L_2|\ge  \frac{\L_-}{2} \left(\frac{m_1}{m_2}\sqrt{(1-\d)\a_-}-1\right)$ and, as in the proof of~\cite[Proosition III.2]{pinzari18}, one can take
  $\varepsilon_0<\frac{\L_-}{2} \left(\frac{m_1}{m_2}\sqrt{(1-\d)\a_-}-1\right)$ in order that the denominators of the functions ${\rm c}_1^*$, ${\rm c}_2$, ${\rm c}^*_2$ in~\cite[Appendix A]{pinzari18} do not vanish, and so small that collisions are excluded.
} have

\begin{proposition}[{\cite[Section III and Appendix A]{pinzari18}}] \label{non res orc}
 One can find $\varepsilon_0>0$, depending only on $\L_-$, $\delta$, $\a_-$, $m_1$, $m_2$ such that
the function $\HH_{\textrm{\sc rps}_{\pi}}$ in~\equ{rps Ham real}
 is real--analytic\footnote{Namely, analytic on a complex neighborhood of ${\cal M}_{\textrm{\sc rps}_\pi, \varepsilon_0}$ and real--valued on ${\cal M}_{\textrm{\sc rps}_\pi, \varepsilon_0}$.} for $(\L,\l,\eta,\xi,p,q)\in{\cal M}_{\textrm{\sc rps}_\pi, \varepsilon_0}$.
 In addition, for any $s\in {\mathbb N}$, there exists a positive number $\a^{\#}$ such that, if $\a_+<\a^{\#}$,
there exists a positive number $\varepsilon_1<\varepsilon_0$ and a real--analytic canonical transformation
\beqano
\phi_{\textrm{\sc bnf}}:\quad
(\L,\ovl\l, \ovl \eta, \ovl\xi, \ovl p,  \ovl q)\in {\cal M}_{\textrm{\sc rps}_\pi, \varepsilon_1}\to (\L, \l, \eta, \xi, p, q)
\in {\cal M}_{\textrm{\sc rps}_\pi, \varepsilon_0}
\eeqano
which carries $(\ovl\eta, \ovl\xi, \ovl p, \ovl q)=0$ to $(\eta, \xi, p, q)=0$ for all $(\ovl\L, \ovl\l)\in {\cal L}\times {\mathbb T}^2$,
such that, if
\beqa{normalized orc}
\HH_{\textrm{\sc bnf}}:=\HH_{\textrm{\sc rps}_{\pi}}\circ\phi_{\textrm{\sc bnf}}=\hh_{\rm k}(\L)+\m f_{\textrm{\sc bnf}}(\L,\ovl\l,\ovl \eta, \ovl \xi, \ovl p, \ovl q)
\eeqa
then
the averaged perturbing function
$$f_{\textrm{\sc bnf}}^{\rm av}(\L,\ovl \eta, \ovl \xi, \ovl p, \ovl q):=\frac{1}{(2\p)^2}\int_{{\mathbb T}^2} f(\L,\ovl\l,\ovl \eta, \ovl \xi, \ovl p, \ovl q)d\ovl\l_1d\ovl\l_2$$
``is in Birkhoff Normal Form of order $s$'', namely:
\beqano
f_{\textrm{\sc bnf}}^{\rm av}
=C_0(\L)+\O\cdot\ovl\t+\frac{1}{2}\ovl\t\cdot{\rm T}(\L)\ovl\t+{\mathbb 1}_{s\ge 3}\sum_{j=3}^{s}\cP_j(\ovl\t;\L)+{\rm O}_{2s+1}(\ovl\eta,\ovl\xi, \ovl p, \ovl q;\L)
\eeqano
where $\O(\L)=(\O_1(\L), \O_2(\L), \O_3(\L))$; $\cP_j(\ovl\t;\L)$ are homogeneous polynomials of degree $j$ in $\ovl\t:=\left(\frac{\ovl\eta_1^2+\ovl\xi_1^2}{2},\ \frac{\ovl\eta_2^2+\ovl\xi_2^2}{2}\ ,\ \frac{\ovl p^2+\ovl q^2}{2}\right)$ and the determinant of the $3\times 3$ matrix ${\rm T}(\L)$ does not identically vanish.
Moreover, $\phi_{\textrm{\sc bnf}}$ leaves $\GG_{\textrm{\sc rps}_{\pi}}$ unvaried, meaning that the function   
\beqano\ovl\GG:=\L_1-\L_2-\frac{\ovl\eta_1^2+\ovl\xi_1^2}{2}+\frac{\ovl\eta_2^2+\ovl\xi_2^2}{2}+\frac{\ovl p^2+\ovl q^2}{2}\eeqano
is still a first integral to $\ovl{\rm H}$.
\vskip.1in
\noindent

\end{proposition}

\subsection{Hyperbolicity}\label{The hyperbolic character}
The hyperbolic character appears using a  set of canonical coordinates,  named {\it perihelia reduction $(${\sc p}--coordinates$)$}. This is a  further set of canonical coordinates 
\beqa{pgeneral}\textrm{{\sc p}}:=\Big(\widehat{\textrm{{\sc p}}}, (\ZZ, \GG, \zeta,  {\rm g})\Big)\in {\mathbb R}^{3n-2}\times {\mathbb T}^{3n-2}\times {\mathbb R}^{2}\times {\mathbb T}^{2}\eeqa
  performing full reduction of SO(3) invariance for a $n$--particle system, which, in addition keeps regular for co--planar motions. The $\textrm{{\sc p}}$--coordinates  have been firstly introduced in~\cite{pinzari18a}, to which we refer for the proof of their canonical character.
We remark that in~\equ{pgeneral},
  $\GG$, $\ZZ$ and $\zeta$ are the same  as in {\sc jrd} in~\equ{coordinates}.  The coordinate ${\rm g}$, conjugated to $\GG$, is not the same as in~\equ{coordinates}, but of course $(\ZZ, \zeta,  {\rm g})$ are again ignorable and $\GG$ is constant in SO(3) invariant systems. 
For the 3--body problem, namely, $n=2$, 
the $8$--plet ${\widehat{\textrm{{\sc p}}}}$ is given by
  $${\widehat{\textrm{{\sc p}}}}:=(\L_1,\L_2,\GG_2,\Theta, \ell_1,\ell_2, {\rm g}_2, \vartheta)$$
with $\Lambda_j$, $\ell_j$, $\GG_2$ as in~\equ{coordinates}.
To define $\Theta$,  ${\rm g}$, $\vartheta$ and ${\rm g}_2$, we assume that 
  \beqa{good nodes}
\n_1:=k^\ppt\times \CC\,,\qquad {\rm n}_1:={\rm C}\times {\rm P}^\ppu,\qquad   \n_2:={\rm P}^{(1)}\times {\rm C}^\ppd, \quad {\rm n}_2={\rm C}^\ppd\times {\rm P}^\ppd
\eeqa
do not vanish. Note that $\n_1$ in~\equ{good nodes} is the same as in~\equ{Dep nodes}. We\footnote{ The second equality in the first equation in~\equ{belle*} is implied by ${\rm C}=\CC^\ppu+{\rm C}^\ppd$ and $\CC^\ppu\cdot{\rm P}^{(1)}=0 $.} let (under the same notations as in the previous section)
  \beqa{belle*}
\begin{array}{llllrrr}
\begin{array}{lrrr}
 \Theta:={\rm C}\cdot {\rm P}^{(1)}={\rm C}^\ppd\cdot {\rm P}^{(1)} \\
\end{array}
\qquad
\left\{
\begin{array}{lrrr}
 \vartheta:=\a_{{\rm P}^{(1)}}({\rm n}_1, \n_2)&\\
 {\rm g}:=\a_{{\rm C}}(\n_1, {\rm n}_1)\qquad&\\
 {\rm g}_2:=\a_{{\rm C}^\ppd}(\n_2, {\rm n}_2)&\\
\end{array}
\right.
\end{array}
\eeqa
We now describe the r\^ole of the {\sc p}--coordinates in the Hamiltonian~\equ{helio}. We denote as

\beqano{\rm H}_{{\textrm{{\sc p}}}}={\rm h}_{{\rm  k}}(\L_1,\L_2)+\m f_{\textrm{\sc p}}(\L_1,\L_2,\GG_2,\Theta;\ell_1,\ell_2,{\rm g}_2,\vartheta; \GG)\eeqano  where
  \beqano
{\rm h}_{{\rm  k}}(\L_1,\L_2)=-\frac{{\rm  m}_1^3{\rm  M}_1^2}{2\L_1^2}-\frac{{\rm  m}_2^3{\rm  M}_2^2}{2\L_2^2},\qquad f_{\textrm{\sc p}}=-\frac{m_1m_2}{|x_{\textrm{\sc p}}^\ppu-x_{\textrm{\sc p}}^\ppd|}+\frac{y_{\textrm{\sc p}}^\ppu\cdot y_{\textrm{\sc p}}^\ppd}{m_0}.
\eeqano
the Hamiltonian~\equ{helio} expressed in terms of ${\textrm{\sc p}}$, and
\beqano f^{\rm av}_{\textrm{\sc p}}:=\frac{1}{(2\p)^2}\int_{[0,2\p]^2}f_{\textrm{\sc p}}d\ell_1d\ell_2\eeqano
the doubly averaged perturbing function. We look at the expansion

\beqano{f^{\rm av}_{\textrm{\sc p}}}=-\frac{m_1 m_2}{a_2}\Big(1+\alpha^2\PP+{\rm O}(\alpha^3)\Big)\eeqano
where $\alpha:=\frac{a_1}{a_2}$ is the semi--major axes ratio. We focus on the function $\PP$. Let ${\cal L}$ as in~\equ{L0};  $c\in(0,1)$, and put

\beqa{LpGp}
\cL_{\textrm{\sc p}}(\GG)&:=&\Big\{\L=(\L_1,\L_2)\in {\cal L}:\quad \L_1>\GG+\frac{2}{ c}\sqrt{\a_+}\L_2\nonumber\\
&& \hspace*{1em} 5\L_1^2\GG -(\GG+\frac{2}{ c}\sqrt{\a_+}\L_1)^2 (4 \GG+\frac{2}{ c}\sqrt{\a_+}\L_1)>0,\nonumber\\
&& \hspace*{1em}5\Lambda_1^2\GG-(\GG+\Lambda_2)(4\GG+\Lambda_2)>0\nonumber\\
&& \hspace*{1em} \L_2>\GG\ ,\quad \L_1> 2\GG\Big\}\\\label{Lp}
\cG_{\textrm{\sc p}}(\L_1,\L_2,\GG)&:=&\Big\{\GG_2:\ \max\{\frac{2}{ c}\sqrt{\a_+}\L_2, \GG\}<\GG_2<\Lambda_2
\Big\}
\\
\cB_{\textrm{\sc p}}(\GG)&:=&\Big\{(\Theta,\vartheta): \ |\Theta|< \frac{\GG}{2}\,,\  |\vartheta|< \frac{\p}{2}\Big\}\nonumber
\eeqa
and finally
$$\cA_{\textrm{\sc p}}(\GG):=\Big\{(\L_1, \L_2, \GG_2):\ (\L_1, \L_2)\in \cL_{\textrm{\sc p}}(\GG)\,,\  \GG_2\in {\cal G}_{\textrm{\sc p}}(\L_1, \L_2)\Big\}$$
Moreover, we let
\beqa{equilibria}
\begin{array}{llll}
{\cal N}(\GG):=\cA_{\textrm{\sc p}}(\GG)\times{\mathbb T}^3\times \cB_{\textrm{\sc p}}(\GG)\,,\quad {\cal N}_{0}(\GG):= \cA_{\textrm{\sc p}}(\GG)\times{\mathbb T}^3\times\big\{0\,,0\big\}.
\end{array}
\eeqa
Note that phase points in ${\cal N}_{0}$ has the geometrical meaning of co--planar motions with the outer planet in retrograde motion.
\begin{proposition}[{\cite[Section IV]{pinzari18}}]\label{hyperHAM}
The 4 degrees of freedom Hamiltonian ${\rm H}_{\textrm{\sc p}}$
is real--analytic in 
$
{\cal N}
$.
It has an equilibrium on ${\cal N}_{0}$. Such equilibrium turns to be hyperbolic\footnote{In~\cite{pinzari18a} a slightly more general result is proved: the equilibrium is hyperbolic when ${\cal L}_{\textrm{\sc p}}$ in~\equ{LpGp} is defined without the inequality \beqa{L2cond}5\Lambda_1^2\GG-(\GG+\Lambda_2)(4\GG+\Lambda_2)>0\eeqa
and
${\cal G}_{\textrm{\sc p}}$ in~\equ{Lp} is taken to be $\Big\{\GG_2:\ \max\{\frac{2}{ c}\sqrt{\a_+}\L_2, \GG\}<\GG_2<\min\{ \L_2, \GG^{\star}\}\Big\}$, with
$\GG^\star$ the unique root of the polynomial
$\GG_2\to 5\Lambda_1^2\GG-(\GG+\GG_2)(4\GG+\GG_2)$.
But as~\equ{L2cond} ensures $\Lambda_2<\GG^\star$, under such restriction,  ${\cal G}_{\textrm{\sc p}}$ can be taken as in~\equ{Lp}.
 } for $\PP$.

\end{proposition}

\subsection{Existence and co--existence of two families of tori}\label{Existence and co--existence of two families of tori}

Theorem~\ref{stable toriREF}   and Theorem~\ref{thm:simplifiedHYPER} can now be used to prove the existence of both full--dimensional and whiskered, co--dimension 2 tori in the three--body problem. Indeed:

\begin{itemize}
\item[--] 
Under conditions~\equ{epsmuREF}, by Theorem~\ref{stable toriREF}, an invariant\footnote{More precisely, Theorem~\ref{stable toriREF} is applied to the Hamiltonian $\HH_{\textrm{\sc bnf}}$ in~\equ{normalized orc}, hence with
$$n_1=2\,,\ n_2=3\,,\ V={\cal L}\,,\ \varepsilon=\varepsilon_1\,,\ H_0=\hh_{\rm k}\,,\ P=f_{\textrm{\sc bnf}}$$
corresponding to the image under $\overline\phi$ of the invariant set obtained through the thesis of Theorem~\ref{stable toriREF}.
} set ${\cal F}\subset{\cal M}_{\varepsilon}
$ for the Hamiltonian $\HH_{\textrm{\sc rps}_{\pi}}$  with $5$--dimensional frequencies is found, whose measure satisfies
 \beqa{measstableREFNEW}
\meas {\cal M}_{\varepsilon}>\meas {\cal F}> \Big(1- C_* {\varepsilon}^{\frac{1}{2}+\ovl s}\Big) \meas {\cal M}_{\varepsilon}
\eeqa
where $\ovl s=s-2$.
\item[--] Under conditions~\equ{epsmuHYPER} with $a=\alpha_+$, by Theorem~\ref{thm:simplifiedHYPER}, for any $\GG\in {\mathbb R}_+$, one finds an invariant set ${\cal H}(\GG)\subset  {\cal N}_{0}(\GG)$  with $3$--dimensional frequencies for
$\HH_{\textrm{\sc p}}$ and equipped with $4$--dimensional stable and unstable manifolds\footnote{Theorem~\ref{thm:simplifiedHYPER} is applied to the Hamiltonian $\HH_{\textrm{\sc p}}$ of Proposition~\ref{hyperHAM}, hence with
$$n_1=2\,,\ n_2=1\,,\ V=\cA_{\textrm{\sc p}}(\GG)$$
$$H_0=\hh_{\rm k}-\frac{m_1 m_2}{a_2}\,,\ P_0=-\frac{m_1 m_2}{a_2}\alpha^2\PP\,,\ P_1=-\frac{m_1 m_2}{a_2}{\rm O}(\alpha^3)\,,\ a=\alpha_+$$},
  whose measure satisfies
\beqa{meas est*HYPERNEW}
\meas  {\cal N}_{0}(\GG)>\meas{\cal H}(\GG)> \Big(1- C_* \sqrt{ \alpha_+ }\Big) \meas  {\cal N}_{0}(\GG)\ .
\eeqa
\end{itemize}
In the next, we show that the invariant sets ${\cal F}$ and ${\cal H}(\GG)$ constructed above ``have a common domain of existence''. We have to make this assertion more precise, mainly because 
 ${\cal F}$ and ${\cal H}(\GG)$ have been constructed with different formalisms.

\vskip.1in
\noindent
Let \beqa{change}\phi_{\textrm{\sc rps}_\p}^{\textrm{\sc p}}:\ {\textrm{\sc rps}_\p}\to \textrm{\sc p}\eeqa the canonical change of coordinates between $\textrm{\sc rps}_\p$ and $\textrm{\sc p}$, well defined in a full measure set.

\noindent
Let ${\mathbb G}_*$, ${\mathbb G}_0$ the respective images  under the function~\equ{Grps}:
$${\mathbb G}_0:=\GG_{\textrm{\sc rps}_\p}\left({\cal M}_{\varepsilon}\right)\,,\qquad {\mathbb G}_*:=\GG_{\textrm{\sc rps}_\p}\left({\cal F}\right)$$
of the sets ${\cal M}_{\varepsilon}$, ${\cal F}$. As ${\cal F}\subset {\cal M}_{\varepsilon}$,  then 
${\mathbb G}_*\subset{\mathbb G}_0$.
For any $\GG_0\in {\mathbb G}_0$, $\GG_*\in {\mathbb G}_*$, let 
$${\cal M}_{\varepsilon}(\GG_0):={\cal M}_{\varepsilon}\cap \{\GG_{\textrm{\sc rps}_\p}=\GG_0\}\,,\qquad {\cal F}(\GG_*):={\cal F}\cap  \{\GG_{\textrm{\sc rps}_\p}=\GG_*\}$$
${\cal M}_{\varepsilon}(\GG_0)$ and ${\cal F}(\GG_*)$ are invariant sets because $\GG_{\textrm{\sc rps}_\p}$ is conserved along the motions of $\HH_{\textrm{\sc rps}}$. 

\noindent
Define:
$${\cal M}'_{\varepsilon}(\GG_0):=\phi_{\textrm{\sc rps}_\p}^{\textrm{\sc p}}\left({\cal M}_{\varepsilon}(\GG_0)\right)\,,\qquad
{\cal F}'(\GG_*):=\phi_{\textrm{\sc rps}_\p}^{\textrm{\sc p}}\left({\cal F}(\GG_*)\right)\,.
$$
At the cost of eliminating zero--measure sets from ${\mathbb G}_0$, ${\mathbb G}_*$, 
  the sets
${\cal F}'(\GG_*)$, 
${\cal M}'_{\varepsilon}(\GG_0) $
are well--defined, for all $\GG_0\in {\mathbb G}_0$, $\GG_*\in {\mathbb G}_*$.
Then split
$${\cal M}'_{\varepsilon}(\GG_0)= \widehat{\cal M}'_{\varepsilon}(\GG_0)\times \{\GG=\GG_0\,,\ {\rm g}\in {\mathbb T}\}\qquad 
{\cal F}'(\GG_*)=\widehat{\cal F}'(\GG_*)\times \{\GG=\GG_*\,,\ {\rm g}\in {\mathbb T}\}
$$

\noindent
The volume--preserving property of $\phi_{\textrm{\sc rps}_\p}^{\textrm{\sc p}}$ in~\equ{change}, the monotonicity of the Lebesgue integral and the bounds in~\equ{measstableREFNEW} guarantee that 
 \beqa{stablereducedtorimeas}
\meas \widehat{\cal M}'_{\varepsilon}(\GG_*)>\meas \widehat{\cal F}'(\GG_*)> \Big(1- C_1 {\varepsilon}^{\frac{1}{2}+\ovl s}\Big) \meas \widehat{\cal M}'_{\varepsilon}(\GG_*)\qquad \forall\ \GG_*\in{\mathbb G}_*\,.
\eeqa
with some $C_1>0$. 
\vskip.1in
\noindent
Recall now the definition of ${\cal N}(\GG)$, ${\cal N}_{0}(\GG)$ in~\equ{equilibria} and ${\cal H}(\GG)$ in~\equ{meas est*HYPERNEW}. 
 The main result of the paper is the folowing
\begin{theorem}\label{main} 
 Let  $\s>0$ half--integer. There exist $\varepsilon_*$,  $c_0\in (0, 1)$
such that, if  $\varepsilon<\varepsilon_*$, $\GG_*\in {\mathbb G}_*$,  $\GG_*>c_0^{-1}\varepsilon^2$, $\alpha_+\le c_0 \varepsilon^{12}$ and $\mu$ verifies~\equ{epsmuREF},~\equ{epsmuHYPER} with $a=\alpha_+$ and $s=\sigma+\frac{7}{2}$, then
there exists a non--empty set 
${\cal A}_{\star}(\GG_*)
$
such that, letting
$${\cal Q}(\GG_*):={\cal A}_{\star}(\GG_*)\times{\mathbb T}^3\times \cB_1(\varepsilon, \GG_\star)\,,\qquad {\cal Q}_0(\GG_*):={\cal A}_{\star}(\GG_*)\times{\mathbb T}^3\times \{(0, 0)\}
$$ and denoting 
$\widehat{\cal F}_*'(\GG_*)$, $\widehat{\cal H}_*(\GG_*)$
the respective intersections of $\widehat{\cal F}'(\GG_*)$, $\widehat{\cal H}(\GG_*)$
with ${\cal Q}(\GG_*)
$, ${\cal Q}_0(\GG_*)
$
then $\widehat{\cal F}_*'(\GG_*)$, $\widehat{\cal H}_*(\GG_*)$ are non--empty and in fact verify
\beqa{first bound}&&\meas{\cal Q}(\GG_*)\ge\meas\widehat{\cal F}_*'(\GG_*)\ge \left(1-\frac{\varepsilon^\s}{\phantom{c_0}\varepsilon_*^{\s\phantom{2}}}\right) \meas{\cal Q}(\GG_*)\\
\label{second bound}&&\meas{\cal Q}_{0}(\GG_*)\ge \meas \widehat{\cal H}_*(\GG_*)\ge \left(1-\frac{\alpha_+}{c_0\varepsilon^{12}}\right)\meas{\cal Q}_{0}(\GG_*)\,.\eeqa

\end{theorem}
The proof of Theorem~\ref{main}  relies on some technical result (Propositions~\ref{inclusion},~\ref{inclusion1} and~\ref{lem: measure} below) which we now state and prove later.

\begin{proposition}\label{inclusion}
Let, for a suitable pure number $\underline k\in (1, 2)$, $\Lambda_-<\GG$, $k_-\le \underline k$  $k_+\ge 2 $, $\alpha_+\le \frac{c^2}{16}$. Choose
$\L_+$ as the unique value of $\L_2> \GG$ such that $\cC$ and the straight line $\L_1=2 \L_2$ meet at $(\L_1, \L_2)=(2\L_+, \L_+)$. 
Let
$${\cal L}_0(\GG):=\Bigg\{(\L_1, \L_2):\quad \GG\le \L_2\le \L_+\ ,\quad  (\GG+\L_2)\sqrt{\frac{4\GG +\L_2}{5\GG} }<\L_1< \min\{k_+\,\L_2\,,\ 2\L_+\}\Bigg\}$$
$$\cA_{0}(\GG):=\Big\{(\L_1, \L_2, \GG_2):\ (\L_1, \L_2)\in \cL_{0}(\GG)\,,\  \GG_2\in {\cal G}_{\textrm{\sc p}}(\L_1, \L_2)\Big\}$$
Then
the set \beqano
\begin{array}{llll}
{\cal N}_{0}(\GG):=\cA_{0}(\GG)\times{\mathbb T}^3\times \cB_{\textrm{\sc p}}(\GG)
\end{array}
\eeqano
is a subset of ${\cal N}(\GG)$.\end{proposition}

\begin{proposition}\label{inclusion1} There exists $c_1\in (0, 1)$ depending only on $\L_+/\GG$, $\L_-/\GG$ such that, letting, for any  $\g< c_1^2\varepsilon^2$,
$$\cL_{1}(\GG):=\Big\{(\L_1, \L_2)\in {\cal L}\ ,\ |\L_1-\L_2-\GG|<c^2_1\varepsilon^2\Big\}$$
$$ \cG_{1}(\L_2):=\Big\{\GG_2:\ \L_2-c^2_1\varepsilon^2<\GG_2<\L_2-\g
\Big\}$$
$$\cA_{1}(\GG):=\Big\{(\L_1, \L_2, \GG_2):\quad (\L_1, \L_2)\in \cL_{1}\ ,\quad \GG_2\in \cG_{1}(\L_2)\Big\}$$
$$\cB_{1}(\GG, \varepsilon):=\Big\{(\Theta, \vartheta):\quad \Theta^2<c^2_1\GG\varepsilon^2\,,\ \vartheta^2<c^2_1\frac{\varepsilon^2}{\GG}\Big\}$$
then
the set \beqano
\begin{array}{llll}
{\cal N}_{1}(\GG, \varepsilon):=\cA_{1}(\GG)\times{\mathbb T}^3\times \cB_{1}(\GG, \varepsilon)
\end{array}
\eeqano
is a subset of $ \widehat{\cal M}'_{\varepsilon}(\GG)$.
\end{proposition}

\begin{proposition}\label{lem: measure} Assume $\GG\ge 10c_1^2\varepsilon^2$ and ${\a_+}<\frac{c^2}{16}$. Then
${\cal A}_0(\GG)$ and ${\cal A}_1(\GG)$ have a non--empty intersection ${\cal A}_\star(\GG)$, verifying
$$\meas (\cA_\star(\GG))\ge \frac{9}{10}(c_1^2\varepsilon^2-\g)c_1^4\varepsilon^4$$ 

\end{proposition}

\noindent
We prove how Theorem~\ref{main}  follows from the above propositions.  
 ${\cal Q}(\GG_*)$ is a subset of $\widehat{\cal M}'_{\varepsilon}(\GG_*)$
and
${\cal N}_{\varepsilon}(\GG_*)$, and
$$\meas {\cal Q}(\GG_*)= C_1 \varepsilon^8 =C_2\varepsilon^2 \meas \widehat{\cal M}'_{\varepsilon}\,.$$
The bound in~\equ{stablereducedtorimeas} guarantees that 
\beqano
\meas \left(\widehat{\cal M}'_{\varepsilon}(\GG_*)\setminus \widehat{\cal F}'(\GG_*)\right)<C_3 {\varepsilon}^{\frac{1}{2}+\ovl s} \meas \widehat{\cal M}'_{\varepsilon}(\GG_*)\qquad \forall\ \GG_*\in{\mathbb G}_*\,.
\eeqano
On the other hand, if $\widehat{\cal F}'_{\varepsilon}(\GG_*)\cap{\cal Q}(\GG_*) $ was empty, we would have
\beqano
\meas \left(\widehat{\cal M}'_{\varepsilon}(\GG_*)\setminus \widehat{\cal F}'(\GG_*)\right)\ge 
\meas {\cal Q}(\GG_*)= C_2\varepsilon^2 \meas \widehat{\cal M}'_{\varepsilon}(\GG_*)
\eeqano
which contradicts the previous inequality if $\ovl s> \frac{3}{2}$ and $\varepsilon$ is small.
Finally, if $\GG_*\in{\mathbb G}_*$,
\beqano
\meas \left({\cal Q}_{\varepsilon}(\GG_*)\setminus \widehat{\cal F}'(\GG_*)\right)\le \meas \left(\widehat{\cal M}'_{\varepsilon}(\GG_*)\setminus \widehat{\cal F}'(\GG_*)\right)<C_3 {\varepsilon}^{\frac{1}{2}+\ovl s} \meas \widehat{\cal M}'_{\varepsilon}(\GG_*)=C_4 {\varepsilon}^{\ovl s-\frac{3}{2}} \meas{\cal Q}_{\varepsilon}(\GG_*)
\eeqano
and we have~\equ{first bound} with $\s=\ovl s-\frac{3}{2}=s-\frac{7}{2}$, with $s\ge 4$. The proof of~\equ{second bound} is similar.
\vskip.1in
\noindent

\subsection{Proof of Propositions~\ref{inclusion},~\ref{inclusion1} and~\ref{lem: measure}}\label{Proof of Propositions}

\proof\,{\bf of Proposition~\ref{inclusion}} We only need to prove that
${\cal L}_0(\GG)\subset {\cal L}_{\textrm{\sc p}}(\GG)$.
We switch to the coordinates

$$y:=\frac{\L_1}{\GG}\,,\qquad x:=\frac{\L_2}{\GG}\,.$$
We denote as
 ${\cal X}_{\textrm{\sc p}}:=\GG^{-1}\cL_{\textrm{\sc p}}$ the domain of $(y, x)$, and as
$$x_-:=\frac{\Lambda_-}{\GG}\,,\ \quad x_+:=\frac{\Lambda_+}{\GG}
$$
 ${\cal X}_{\textrm{\sc p}}$ can be written as the intersection of the three sets:

 \beqano
{\cal X}_1&:=&\Bigg\{(y, x):\quad 1\le x\le x_+\ ,\quad y>2\,,\  \max\{k_-\,x, (1+x)\sqrt{\frac{4 +x}{5} }\}<y< k_+\,x\Bigg\}\nonumber\\
{\cal X}_2&:=&\Bigg\{(y, x):\quad 1\le x\le x_+
\ ,\quad y>1+\frac{2}{ c}\sqrt{\a_+}x \Bigg\}\nonumber\\
{\cal X}_3&:=&\Bigg\{(y, x):\quad 1\le x\le x_+\,,\quad    y>2\ ,\quad 5y^2 -(1+\frac{2}{ c}\sqrt{\a_+}y)^2 (4 +\frac{2}{ c}\sqrt{\a_+}y)>0\Bigg\}
\eeqano
We prove  ${\cal X}_0:=\GG^{-1}\cL_0$ is a subset of all of them.
The curve
$${\cal C}:\qquad y=(1+x)\sqrt{\frac{4 +x}{5}} \qquad x\ge 1$$
passes through $P_0=(1, 2)$.  We denote as $\underline k$ the slope of the straight line $y=kx$ which is tangent at ${\cal C}$  at $P_0$. The slope of the straight line $y=kx$ through $P_0$ is obviously $\overline k=2$. 
We assume that
$$k_-\le \underline k\,,\qquad k_+\ge \overline k
$$
and choose $(x_+, y_+)$ as the only $(x, y)$ with $x>1$ such that $\cC$ meets $y=2 x$ at $(x, y)$.
Under such assumptions, we have:
$${\cal X}_1=\Bigg\{(y, x):\quad 1\le x\le x_+\ ,\quad  (1+x)\sqrt{\frac{4 +x}{5} }<y< k_+\,x\Bigg\}\supset{\cal X}_0$$
The straight line which is tangent at $\cC$ at $P_0=(1, 2)$ has equation
$$ y=\frac{6}{5}x+\frac{4}{5}$$
Since we  $\a_+<\frac{{c^2}}{4}$, $x>1$ and $\cC$ is convex, we have
$$1+\frac{2}{ c}\sqrt{\a_+}x\le 1+x
\le \frac{6}{5}x+\frac{4}{5}\le (1+x)\sqrt{\frac{4 +x}{5} }$$
This shows that ${\cal X}_2\supset {\cal X}_0$.
As for ${\cal X}_3$, we note that
for
$$\alpha_+\le \frac{c^2}{16}$$
it is
$$5y^2 -(1+\frac{2}{ c}\sqrt{\a_+}y)^2 (4 +\frac{2}{ c}\sqrt{\a_+}y)\ge 5y^2 -\left(1+\frac{y}{ 2}\right)^2 \left(4 +\frac{y}{ 2}\right)=\frac{1}{4}(y-2)(y-y_-)(y_+-y)\,.$$
with
$$y_\pm=13\pm \sqrt{185}\,.$$
As $y_-<0$ and  $(y-2)(y_+-y)\ge 0$  on ${\cal X}_0$,  we have that ${\cal X}_3\supset {\cal X}_0$. $\qquad \square$

\begin{figure}
\centering{
\includegraphics{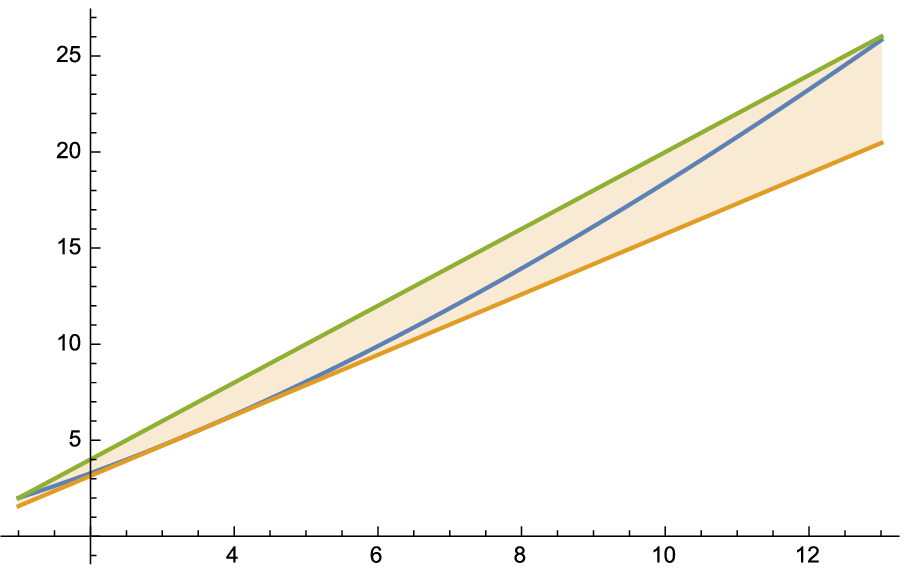}}
\caption{The blue curve is ${\cal C}$; the orange line has slope $k_-$, the green one has slope $k_+$ (\textsc{Mathematica}).\label{coexistence1}}
\end{figure}
\begin{figure}
\centering{
\includegraphics{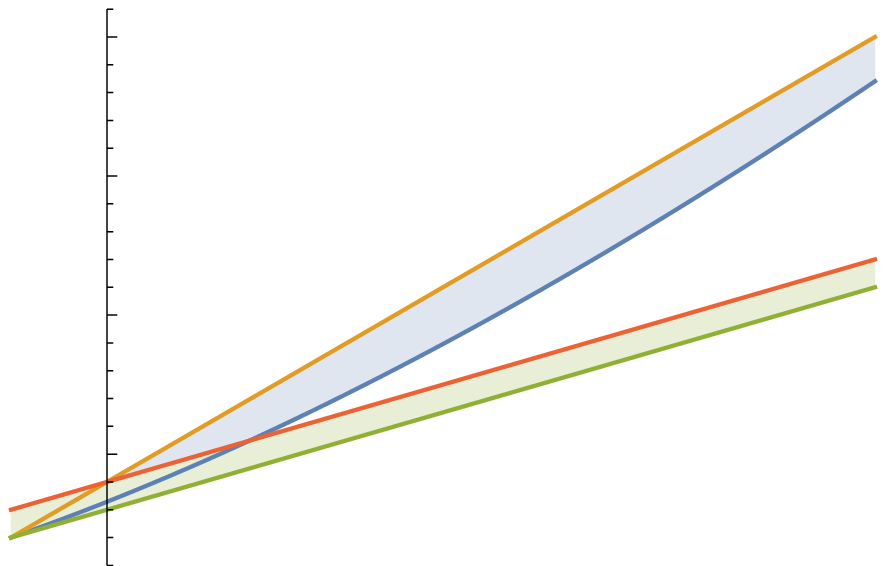}}
\caption{
$\cL_0(\GG)$ (blue) and $\cL_1(\GG)$ (green).\label{coexistence}}
\end{figure}

\begin{remark}\rm The numbers $\underline k$, $\L_+$ of Proposition~\ref{inclusion} can be chosen as
$$\L_+=\frac{\GG}{2}\left(13+ \sqrt{185}\right)\,,\quad \underline k=\frac{1}{4}\sqrt{\frac{3}{10}(69+11\sqrt{33})}\sim 1.57$$
\end{remark}
$\L_+$ is related to the number $x_+$  computed along the proof via $\L_+=x_+\GG$. $\underline k$ is defined as the slope of the straight line $y= kx $ which is tangent at $\cC$. We can compute it eliminating $y$ between the two equations, we obtain the  cubic equation
\beq{1st}x^3+(6-5   k^2)x^2+9x+4=0.\eeq
The tangency condition is imposed identifying this equation with
 \beq{2nd}(x-a)^2(x-b)=0
\eeq
where $a$ is the abscissa of the tangency point. 
Equating the respective coefficients of~\equ{1st} and~\equ{2nd}, we obtain
\beqa{syst}\arr{ -(b+2a)=6-5  k^2\\
  2ab+ a^2=9\\
  -a^2b=4
 }\eeqa
Eliminating $b$ through the second and the third equations, we obtain
$$a^3-9 a-8=0$$
which has the following three roots: 
$$a_0=-1,\qquad a_{\pm}=\frac{1\pm\sqrt{33}}{2}.$$
The only admissible  value is then
$$a=a_+=\frac{1+\sqrt{33}}{2}\ .$$
In correspondence of this value for $a$, solving  the system in~\equ{syst}, we find
$$b=\frac{-17+\sqrt{33}}{32},\qquad k=\frac{1}{4}\sqrt{\frac{3}{10}(69+11\sqrt{33})}=\underline k\,.\quad \square $$

\proof\,{\bf of Proposition~\ref{inclusion1}} From~\equ{PR}, we get
$$|z|^2=\eta_1^2+\xi_1^2+\eta_2^2+\xi_2^2+p^2+q^2=2(\GG+\GG_2-\GG_1)+2(\L_1-\GG_1)+2(\L_2-\GG_2)\,.$$
From the equality
\beqano
\GG_1&=&\sqrt{\GG^2+\GG_2^2-2\Theta^2+2\sqrt{\GG^2-\Theta^2}\sqrt{\GG_2^2-\Theta^2}\cos\vartheta}\nonumber\\
&=&\GG+\GG_2+{\rm O}\left(
\frac{\Theta^2}{\GG+\GG_2}\right)+{\rm O}\left(
\frac{\Theta^2\GG_2}{\GG(\GG+\GG_2)}\right)+{\rm O}\left(
\frac{\Theta^2\GG}{\GG_2(\GG+\GG_2)}\right)+{\rm O}\left(
\frac{\vartheta^2\GG\GG_2}{\GG+\GG_2}
\right)
\eeqano
and the definition of ${\cal N}_1(\GG)$, the assertion trivially follows.
$\quad \square$

\proof\  {\bf of Proposition~\ref{lem: measure}}   Let $\L_2^\star$ be the abscissa, in the plane $(\Lambda_2, \Lambda_1)$, of the intersection point between the curves
$$ \L_1=(\GG+\L_2)\sqrt{\frac{4 \GG+\L_2}{5\GG}}\ ,\qquad {\rm and}\qquad  \L_1=\L_2+\GG+c_1^2\varepsilon^2\ .$$
Using the coordinate
$x:=\frac{\L_2}{\GG}$. With $x^\star:=\frac{\L_2^\star}{\GG}$, $\theta:=\frac{c_1^2\varepsilon^2}{\GG}$,  $\zeta:=\frac{\g}{\GG}$, where $\zeta<\theta$, the set ${\cal A}_\star(\GG):={\cal A}_0(\GG)\cap {\cal A}_1(\GG)$ has measure 
$$
\meas({\cal A}_\star(\GG))=\GG^3\int_{1+\zeta}^{x^\star}F_1(x) F'_2(x)dx
$$
where
\beqano
F_1(x)&=&\min\Big\{
2x,\  x+1+\theta
\Big\}-(1+x)\sqrt{\frac{4+x}{5}}\nonumber\\
F_2(x)&=&\min\Big\{
\theta-\zeta,\  x-1-\zeta, \ m x-\zeta
\Big\}
\eeqano
and where, for short, we have let $m:=1-\frac{2}{ c}\sqrt{\a_+}$. Then,
\beqa{intHG}
\meas({\cal A}_\star(\GG))\ge\GG^3\int_{1+\zeta}^{x^\star}F_1(x) F_2(x)dx\ .
\eeqa
 To go further, we need a quantitative bound on $x^\star$.  Indeed we have
\begin{claim}\label{lem: zeroes} If $0<\theta<\frac{1}{10}$, then $1+4\theta<x^\star<1+6\theta$.
\end{claim}
The proof of the claim is postponed below, in order not to interrupt the main proof.

\vskip.1in
\noindent
Since we have assumed $\GG\ge 10c_1^2\varepsilon^2$ and 
$\alpha_+\le \frac{c^2}{16}$, 
 then
$\GG\ge \frac{\frac{12}{ c}\sqrt{\a_+}}{1-\frac{2}{ c}\sqrt{\a_+}}c_1^2\varepsilon^2$. In the new variables, this is  $\theta\le \frac{m}{6(1-m)}$. But then $$x^*<1+6\theta\le \frac{1}{1-m}\qquad \Longrightarrow\quad x-1-\zeta\le mx-\zeta\quad \forall\ x<x^*$$ whence 
$$F_2(x)=\left\{
\begin{array}{lll}x-1-\zeta\quad&{\rm if}\quad &1+\zeta\le x\le 1+\theta\\\\
\theta-\zeta &{\rm if} &1+\theta<x\le x^\star
\end{array}
\right.$$
Observe that the second inequality is well put, because $x^\star>1+4\theta$,  as said. The function $F_1(x)$ splits in the same intervals:
\beqa{F2}F_1(x)=\left\{
\begin{array}{lll}2x-
(1+x)\sqrt{\frac{4+x}{5}}
\quad&{\rm if}\quad &1+\zeta\le x\le 1+\theta\\\\
x+1+\theta-(1+x)\sqrt{\frac{4+x}{5}} &{\rm if} &1+\theta<x\le x^\star
\end{array}
\right.\eeqa
Since $\zeta<\theta$, a  lower bound to the integral in~\equ{intHG} is given by
\beqano\int_{1+\zeta}^{x^\star}F_1(x) F_2(x)dx\ge \int_{1+\theta}^{x^\star}F_1(x) F_2(x)dx=(\theta-\zeta) \int_{1+\theta}^{x^\star} F(x)dx\eeqano
with \beqa{zeroes}F(x):=x+1+\theta-(1+x)\sqrt{\frac{x+4}{5}}
\eeqa
 the function in the second line in~\equ{F2}.
Since $F$ is the difference of a linear function and a convex one, 
it is concave. Then we have
$$F(x)\ge F(1)+\frac{F(x^\star)-F(1)}{x^\star-1}(x-1)\qquad \forall\ 1\le  x\le x^\star $$
since $F(x^\star)=0$ and $F(1)=\theta$, this inequality becomes
$$F(x)\ge\frac{x^\star-x}{x^\star-1} \theta\qquad \forall\ 1\le  x\le x^\star $$
hence
$$\int_{1+\theta}^{x^\star}F(x)dx\ge \frac{\theta}{x^\star-1}\int_{1+\theta}^{x^\star}(x^\star-x)dx=\frac{\theta}{2} \frac{(x^\star-1-\theta)^2}{x^\star-1}\ge \frac{9}{10}\theta^2 $$
having used $1+4\theta<x^\star<1+6\theta$.\\
 It remains to prove Claim~\ref{lem: zeroes}.
$x^\star$ is defined as the zero of the function $F$ in~\equ{zeroes}
in the range $(1, +\infty)$.
Multiplying the left hand side of Equation $$x+1+\theta-(1+x)\sqrt{\frac{x+4}{5}}=0$$ 
by $x+1+\theta+(1+x)\sqrt{\frac{x+4}{5}}$, 
we obtain the  algebraic  equation of degree three
$$x^3+x^2+(1+10 \theta) x-1-10\theta-5\theta^2=0$$
which, for $x\ge -1$ is completely equivalent to the initial equation. We aim to apply a bisection argument to the function at left hand side, which we denote as $G(x)$. We have
$$G(1+4\theta)=\theta(
64\theta^2+19\theta-4
)\ ,\qquad G(1+6\theta)=\theta(
216\theta^2+79\theta+4
)$$
and it is immediate to check that
$$G(1+4\theta)<0\qquad G(1+6\theta)>0\qquad \forall\ 0<\theta< \frac{-19+\sqrt{1385}}{
128
}= 0.142\ldots$$
To prove uniqueness, just observe that the function $x\in(0, +\infty)\to G(x)$ is increasing for all $\theta>0$.  This completes the proof.
 $\qquad \square$

\section{Quantitative {\sc kam} theory}\label{KAM theory}
 \subsection{Proof of Theorem~\ref{stable toriREF}}\label{App: stable toriREF}
 The proof of  Theorem~\ref{stable toriREF} is based on an application of~\cite[Proposition 3]{chierchiaPi10}. The method is
is completely analogous to the one used in the proof of~\cite[Theorem 1.3]{chierchiaPi10}, so we shall only say what to change in the proof of~\cite[Theorem 1.3]{chierchiaPi10} in order to obtain the proof of
Theorem~\ref{stable toriREF}.
 The polynomial $N(I,r)$ in the first non numbered formula in
\cite[Section 4]{chierchiaPi10}
 is to be changed as 
 \beq{new N}N(I,r)=P_0(I) +\sum_{i=1}^{m}\O_i(I)r_i+\frac{1}{2}\sum_{i,j=1}^{m} \b_{ij}(I)r_i r_j+{\mathbb 1}_{s\ge 3}\sum_{j=3}^{s}\cP_j(r;I)\,. \eeq 
Equations (60) and (61) in~\cite{chierchiaPi10}
can be modified, respectively, as
\beqa{ovl gam*}
&&\sup_{B_\varepsilon^{2m}\times V_{\rho_0}}|\tilde P_{\rm av}|\le C\varepsilon^{2s+1}\qquad \forall\ 0<\varepsilon<\varepsilon_0\nonumber\\
&&\m<\frac{\varepsilon^{2s+2}}{(\log\varepsilon^{-1})^{2\t+1}}\qquad \ovl\g>\Big(\frac{6(2s+1)}{s_0}\Big)^{\t+\frac{1}{2}}\frac{\sqrt\m(\log\varepsilon^{-1})^{\t+\frac{1}{2}}}{\varepsilon^{s+\frac{1}{2}}}\ .
\eeqa
Analogously to~\cite{chierchiaPi10}, one next applies  Lemma A.1  in~\cite{chierchiaPi10}, but modifying the choice of $K$  as 
\beq{ovl K*}K=\frac{6(2s+1)}{s_0}\log\varepsilon^{-1}\eeq and leaving the other quantities unvaried.
 A bound as in Equation (62) in~\cite{chierchiaPi10} is so obtained, with $H_0$ as in~\cite{chierchiaPi10}, $N(\ovl I,\ovl r)$ as in~\equ{new N}, $\m \widetilde P_{\rm av}(\ovl p, \ovl q, \ovl I)=f_{\textrm{\sc bnf}}(I,\ovl p,\ovl q)-N(\ovl I,\ovl r)$ uniformly bounded by $C\m \varepsilon^{2s+1}$, by~\rm (A${}_2$). Due to the choice of $K$ in~\equ{ovl K*} and the one for $\ovl\g$ in~\equ{ovl gam*},  a bound similar to the one in Equation (63) in~\cite{chierchiaPi10} holds, with the right hand side replaced by $\ovl C\m \varepsilon^{2s+1}$. At this point, one follows the indications 
  in Step 2 of the proof of Theorem 1.3 in
\cite{chierchiaPi10}. Namely, one has to repeat 
the procedure in Steps 5 and 6 of  the proof Theorem 1.4 (previously proved in~\cite{chierchiaPi10}), with the following modification. The annulus $\cA(\varepsilon)$ in Equation (47) in~\cite{chierchiaPi10}
 is to be taken as
 \beqano
 \cA(\varepsilon)=\Big\{J\in {\mathbb R}^{m}:\ \check c_1\varepsilon^{s+\frac{1}{2}}<J_i<\check c_2\varepsilon^2\ ,\quad 1<i<m\Big\}
 \eeqano
and the number $\breve \r$ in 
 Equation (48) in~\cite{chierchiaPi10} is to be replaced with $\breve\r:=\min\{\check c_1\varepsilon^{s+\frac{1}{2}}/2,\ \ovl \r/ {48}\}$. The other quantities remain unvaried. In the remaining  Steps 5 and 6 
 of the proof of Theorem 1.4
in~\cite{chierchiaPi10}  replace the number ``5'' appearing in all the  formulae  with $(2s+1)$ and  $\varepsilon^{n_2/2}$ in Equation (56) (and the formulae below) in~\cite{chierchiaPi10}  with $\varepsilon^{m(s-\frac{3}{2})}$.  $\quad \square$
\vskip.1in
\noindent

 \subsection{Proof of Theorem~\ref{thm:simplifiedHYPER}}
The proof of Theorem~\ref{thm:simplifiedHYPER} proceeds along the same lines as the proof of Theorem~\ref{stable toriREF}, apart for being based on a generalization (Theorem~\ref{two scales KAM} below) of~\cite[Proposition 3]{chierchiaPi10} which now we state.

\noindent
As in~\cite{chierchiaPi10} ${\cal D}_{\gamma_1, \gamma_2,\t}\subset {\mathbb R}^{n}$ denotes the set of vectors $\o=(\o_1, \o_2)\in{\mathbb R}^{n_1}\times {\mathbb R}^{n_2}$ satisfying for any
 $k=(k_1,k_2)\in{\mathbb Z}^{n_1}\times {\mathbb Z}^{n_2}\setminus\{0\}$, inequality 
 \beqa{dioph2sc}
 |\o_1\cdot k_1+\so\cdot k_2|\geq 
\left\{ \begin{array}{l}
\displaystyle\frac{\g_1}{|\tk|^{\t}}\quad {\rm if}\quad  k_1\neq0\ ;\\ \ \\
\displaystyle\frac{\g_2}{| k_2|^{\t}}\quad {\rm if}\quad  k_1= 0\ ,\quad  k_2\neq 0\ .
\end{array}\right.
\eeqa

\begin{theorem}
\label{two scales KAM}
Let ${n_1}$, ${n_2}\in {\mathbb N}$, ${n}:={n_1}+{n_2}$, $\tau>n$, $\fg\geq\ssg>0$, $0<s\leq \frac{\varepsilon}{\ovl\varepsilon+\varepsilon}$, $\r>0$,   $A:={D}_\r\times B^2_{\ovl\varepsilon+\varepsilon}$, and let 
\beqano
{\rm H}
(I,\psi, p, q)={\rm h}(I, pq)+{\rm f}(I,\psi, p, q)
\eeqano 
be real--analytic on $A\times {\mathbb T}_{\ovl s+s}^{n}$.
 {Let $$I=(I_1, I_1)\,,\quad \varpi(I,pq):=\partial_{ {(I,pq)}} {\rm h}(I,pq)=(\omega_1(I_1, I_2, pq), \omega_2(I_1, I_2, pq), \n(I_1, I_2, pq))$$  
 with $\omega_k(I_1, I_2, pq):=\partial_{I_k} {\rm h}(I_1, I_2, pq)$, and assume that
 the map $I\in D_\r\to \o(I,J)$ is a diffeomorphism of $D_\r$ for all $J=pq$, with $(p,q)\in B^2_{\varepsilon}$,
with non singular Hessian matrix $U(I,J):=\partial_{I}^2{\rm h}(I,J)$. 
 }
Let\footnote{The norms will be specified in the next Section~\ref{whiskered KAM thm}.}
\beqano
M\geq \|\partial\omega\|_A\ ,\  \widehat M\geq \|\partial\omega_1\|_A\,,\  \ovl M\geq \|U^{-1}\|_A\ ,\  
E\geq\|{\rm f}\|_{\r,\ovl s+s}\,,\quad \lambda\le \inf |\Re\nu|_A\ .
\eeqano
Assume, for\footnote{\equ{lambda cond} is a simplifying assumption. It may be relaxed.} simplicity,
\beqa{lambda cond}2\frac{s^{\tau}\gamma_2}{6^{\tau}\lambda}\le 1\,.\eeqa
Define
\beqano
&&  \displaystyle  \widehat c:=2^7(n+1)(24)^\tau\ ,\quad  {\widetilde c:=2^{6}}\\
&&  \displaystyle  K:=\frac{32}{s}\ \log_+{\left(\frac{E M^2\,L}{\gamma_1^2}\right)^{-1}}\quad {\rm where}\quad \log_+ a :=\max\{1,\log{a}\}\\ 
&&  \displaystyle   {\widehat\r:=\min\left\{\frac{\gamma_1}{2MK^{\tau+1}}\ ,\ \frac{\gamma_2}{2\widehat M K^{\tau+1}}\,,\ 
 \r\right\}}\,,\quad \widetilde\rho:=\min\left\{\widehat\rho\,,\ \frac{\varepsilon^2}{s}\right\}\\
 \\
&&  \displaystyle  L:=\max\ \Big\{\ovl M\,,\  M^{-1}\,,\ \widehat M^{-1}\Big\}
\\
&&  \widehat E:=\frac{E L}{\widehat\r\widetilde\rho}\ ,\qquad  {\widetilde E:=\frac{E}{\lambda\varepsilon^2}}\,.
\eeqano
Finally, let   $\ovl M_1$, $\ovl M_2$ upper bounds on the norms of the sub--matrices $n_1\times n$, $n_2\times n$ of $U^{-1}$ of the first $n_1$, last $n_2$ rows\footnote{{\it I.e.},
$\displaystyle \ovl M_i\geq \sup_{D_\r}\|T_i\|\ ,\quad i=1,\ 2\ ,\quad \textrm{if}\quad U^{-1}=\left(\begin{array}{lrr}
T_1\\
T_2
\end{array}
\right)\ .$
}.
Assume the perturbation ${\rm f}$ so small that the following ``KAM conditions'' hold
\beq{KAM cond}
\widehat c\widehat E<1\ ,\quad \widetilde c\widetilde E<1
\eeq
Then, for any  $(\p, \k)\in B^2_{\ovl\varepsilon}$ and any
$\o_*\in\O_*(\p\k):=\o({D}, \p\k)\cap\cD_{\gamma_1, \gamma_2,\tau}$, one can find a unique real--analytic embedding
{\small\beqa{perturbed torus}
\phi_{\o_*}:\ {\mathbb T}^{n}\times \{(\p, \k)\}&\to&\Re({D}_r)\times {\mathbb T}^{{n}}\times B^2_{\ovl\varepsilon+r'}\nonumber\\
(\vartheta, \p, \k)&\to&\Big(v(\vartheta, \p, \k; \omega_*), \vartheta+u(\vartheta, \p, \k; \omega_*), 
\p+w(\vartheta, \p, \k; \omega_*), \k+y(\vartheta, \p, \k; \omega_*)
\Big) \eeqa}
  such that $\cM_{\o_*}:=\phi_{{\o_*}}({\mathbb T}^n\times B^2_{\ovl\varepsilon})$ is a real--analytic $(n+2)$--dimensional 
manifold, on which the  ${\rm H}$--flow is analytically conjugated to \beqa{whiskered tori}(\vartheta, \p, \k)\in
{\mathbb T}^{n}\times B^2_{\ovl\varepsilon}
\to (\vartheta+\o_* t,\ \p\to \p e^{-\n_*(\o_*, \p\k)t}, \ \k\to \k e^{\n_*(\o_*, \p\k)t})\,.
\eeqa
In particular, the manifolds $${\rm T}_{{\o_*}}:=\phi_{\o_*}\left({\mathbb T}^n\times\{(0,0)\}\right)$$
are real--analytic $n$--dimensional ${\rm H}$--invariant tori 
embedded in $\Re({D}_r)\times{\mathbb T}^{n}\times B^2_{\ovl\varepsilon}$,
  equipped
with $(n+1)$--dimensional  manifolds
$$\cM_{\rm u}:=\phi_{\o_*}\left({\mathbb T}^n\times\{0\}\times B^1_{\ovl\varepsilon}\right)\ ,\qquad \cM_{\rm s}:=\phi_{\o_*}\left({\mathbb T}^n\times B^1_{\ovl\varepsilon}\times\{0\}\right)$$
on which the motions leave, approach ${\rm T}_{{\o_*}}$ at an exponential rate.
Let $ {\rm T}_{\omega_*, 0}$ denote the projection  of ${\rm T}_{\omega_*}$ on the $(I, \varphi)$--variables, and $\displaystyle  {\rm K}_0:=\bigcup_{\o_*\in\O_*}{\rm T}_{\o_*, 0}$.
Then $ {\rm K}_0$ satisfies the following
 measure\footnote{$\meas_n$ denotes the $n$--dimensional Lebesgue measure.} estimate:
\beq{tori measure}
\meas_{2n}(\Re({D}_r)\times{\mathbb T}^{n}\setminus{\rm K}_0)\leq c_n\Big(\meas({D}\setminus{D}_{\gamma_1, \gamma_2,\tau}\times{\mathbb T}^{n})+\meas(\Re({D}_r)\setminus{D})\times{\mathbb T}^{n}\Big),
\eeq
where  ${D}_{\gamma_1, \gamma_2,\tau}$ denotes the  $\o_0(\cdot, 0)$--preimage of $\cD_{\gamma_1, \gamma_2,\tau}$  and $c_n$ can be taken to be $\displaystyle  c_n=(1+(1+2^8nE)^{2n})^2$.\\
Finally,  the following uniform estimates hold for the embedding $\phi_{\omega_*}$:
\beqa{***}
&&
|  v_1(\vartheta, \p, \k;\o_*)-I_1^0(\p\k; \o_*)|\leq {6} {n}\left(\frac{\ovl M_1}{\ovl M}+\frac{\widehat M}{M}\right)\widehat E\,\widetilde\r\nonumber\\
&& | v_2(\vartheta, \p, \k;\o_*)-I_2^0(\p\k;\o_*)|\leq {6}n\left(\frac{\ovl M_2}{\ovl M}+\frac{\widehat M}{M}\right)\widehat E\,\widetilde\r\ ,\nonumber\\
&& |u(\vartheta, \p, \k;\o_*)|\leq2\,\widehat E\,s\,,\quad  |w(\vartheta, \p, \k;\o_*)|\le 2\,\widehat E\,\varepsilon\nonumber\\
&&  |y(\vartheta, \p, \k;\o_*)|\le 2\,\widehat E\,\varepsilon
\eeqa
where  
$v(\vartheta, \p, \k; \o_*)=(v_1(\vartheta, \p, \k; \o_*), v_2(\vartheta,\p, \k; \o_*))$ and
$I^0(\p\k; \o_*)=(I^0_1(\p\k; \o_*),I^0_2(\p\k; \o_*))\in D$ is the $\o(\cdot, \p\k)$-- pre--image of $\o_*\in\O_*(\p\k)$.
where $r:=8 {n} \widehat E \widetilde\r$, $ r'=2\widehat E\varepsilon$
\end{theorem}
The proof of Theorem~\ref{two scales KAM} is deferred to the next Section~\ref{whiskered KAM thm}. Here we prove how Theorem~\ref{thm:simplifiedHYPER} follows from it.

\noindent
As said, we follow the same ideas  of the proof of Theorem~\ref{App: stable toriREF}, which in turn follows~\cite[Theorem 1.3]{chierchiaPi10}.
 By $(A'_{2})$,
\beqa{init est}P_{\rm av}(I, p, q)=P_0(I, pq)+P_1(I, p, q)\quad {\rm where}\quad |P_1|\le  a \|P_0\|=:\epsilon\,.\eeqa
At this point, proceeding as in~\cite[Proof of Theorem 1.3, Step 1]{chierchiaPi10} but with $\epsilon^5$ replaced by $\epsilon$, under condition
\beqano\m<\frac{\e^{1+\eta}}{(\log(\e^{-1}))^{2\t+1}}\ ,\qquad \overline\g\ge C\Big(\frac{6}{s_0}\Big)^{\t+\frac{1}{2}}\frac{\sqrt{\m}(\log{\e^{-1}})^{\t+\frac{1}{2}}}{\sqrt{\e}}\ ,\eeqano
by an application of~\cite[Lemma~A.1]{chierchiaPi10}, with $\overline K=\frac{6}{s_0}\log{\epsilon^{-1}}$, $r_p=r_q=\e_0$, $r=4\rho=\overline\r:=\min\left\{\frac{\overline\g}{2\overline M\overline K^{\t+1}},\ \r_0\right\}$ (with $\ovl M:=\sup|\partial^2_{I_1}H_0|$), $\r_p=\r_q=\e_0/4$, $\s=s_0/4$, $\ell_1=n_1$, $\ell_2=0$, $m=n_2$ $h=H_0$, $g\equiv 0$, $f=\m P$, $A=\overline D:=\omega_0^{-1}\cD_{\gamma, \tau}$ (where $\omega_0$ is as in $A_1$ and $\cD_{\gamma, \tau}$ is the usual Diophantine set in $\mathbb R^n$, namely the set~\equ{dioph2sc} with $\gamma_1=\gamma_2$), ${B}={B'}=\{0\}$,   $s=s_0$, $\a_1=\a_2=\overline\a=\frac{\overline\g}{2\overline K^\t}$, and $\L=\{0\}$, on the  domain $W_{\ovl v,\ovl s}$  where $\ovl v=(\overline\r/2,\e_0/2)$ and $\ovl s=s_0/2$,
one finds a real--analytic and
 symplectic transformation $\ovl\phi$ which carries $\HH$ to 
 \beqano
\overline H(\overline I,\overline \varphi,\overline p,\overline q)&:=&H\circ\ovl\phi(\overline I,\overline \varphi,\overline p,\overline q)\\
&=&H_0(\overline I)+\m P_{0}(\overline I, \overline p\overline q)+\mu P_1(\overline I,\overline \varphi,\overline p,\overline q)+\widetilde P(\overline I,\overline \varphi,\overline p,\overline q)\nonumber\\
&=&H_0(\overline I)+\m P_{0}(\overline I, \overline p\overline q)+\mu \overline P(\overline I,\overline \varphi,\overline p,\overline q)
\eeqano
where
\beqano\|\widetilde P\|_{\overline v,\overline s}\le{\overline C\m}\max\{\frac{\m \overline K^{2\t+1}}{\overline\g^2},\frac{\m\overline K^\t}{\overline\g}\ e^{-\overline K s_0/2}\}\le {\overline C}\m \epsilon={\overline C}\m a \|P_0\| \ ,\eeqano
whence (by~\equ{init est}) also $\overline P=\m\widetilde P_{\rm av}+\widetilde P$ is bounded by ${C}\m a \|P_0\|$ on  $W_{\overline v,\overline s }$.\\
The next step is to apply Theorem ~\ref{two scales KAM}  to the Hamiltonian $\ovl H$.
Since we can take
\beqano
&&  M=C\ ,\quad 
\sHs=C\m \|P_0\|\ ,\quad \inHs=C(\m \|P_0\|)^{-1}\ ,\quad   E ={C}\m a \|P_0\|\nonumber\\
&&
\bar M_1=C\ , \quad \bar M_2=C(\m \|P_0\|)^{-1}\,,\quad \lambda=C^{-1} \m  \|P_0\|
\eeqano
 the numbers $L$, $K$, $\widehat\rho$ and $\widetilde\rho$ can be bounded, respectively, as
$$L\le C(\m \|P_0\|)^{-1}\ ,\qquad K\leq C\log{(a /\fg^2)^{-1}}$$
and
\beqano
&&\widehat\r\geq c\,\min\Big\{\frac{\fg}{(\log{(a/\fg^2)^{-1}})^{\t+1}}\ ,\ \frac{\sg}{(\log{(a/\fg^2)^{-1}})^{\t+1}} \ ,  \frac{\bar\g}{(\log{\e^{-1}})^{\bar\t+1}}\ ,
\ \r_0\Big\}\nonumber\\
&&\widetilde\r\geq c\,\min\Big\{\frac{\fg}{(\log{(a/\fg^2)^{-1}})^{\t+1}}\ ,\ \frac{\sg}{(\log{(a/\fg^2)^{-1}})^{\t+1}} \ ,  \frac{\bar\g}{(\log{\e^{-1}})^{\bar\t+1}}\ ,
\ \r_0\,,\ \varepsilon^2\Big\}
\eeqano
having let $\gamma_2:=\mu\|P_0\| \ovl\gamma_2$. Condition~\equ{lambda cond}
is trivially satisfied for any $\ovl\gamma<1$, $s\le 6$, while, from the bounds \beqano
\widehat c\widehat E&\le& 
C a \max\left\{
\frac{(\log{(a/\fg^2)^{-1}})^{2(\t+1)}}{\gamma_1^2}\,,\ \frac{(\log{(a/\fg^2)^{-1}})^{2(\t+1)}}{\ovl\gamma_2^2}\,,\ \frac{(\log{\epsilon)^{-1}})^{2(\t+1)}}{\ovl\gamma^2}\,,\ \frac{1}{\rho_0^2}\,,\ \frac{1}{\varepsilon^4} 
\right\}\ ,\ 
\nonumber\\
\widetilde c\widetilde E&\le& C\frac{a}{\varepsilon^2}
\eeqano
one sees that conditions~\equ{KAM cond} hold taking
$$\ovl\gamma=\gamma_1=\ovl\gamma_2=\widehat C\sqrt a\,,\quad a<\widehat C^{-1}\varepsilon^4$$
with a suitable $\widehat C>1$. By the thesis of Theorem  ~\ref{two scales KAM}, we
can find  a set of $n$--dimensional invariant tori ${\cal K}\subset \cP$
whose projection $\cK_0$ on $\cP_0$ satisfies the measure estimate
\beqano
\meas \cP_0\ge \meas\cK_0\ge (1- C'(\overline\g+\g_1+{\g_2}))\meas \cP_0\ge (1- C\sqrt a)\meas \cP_0\ .\quad \square
\eeqano
\subsection{Proof of Theorem~\ref{two scales KAM}}\label{whiskered KAM thm}

We fix the following notations.
\begin{itemize}  
\item[$\bullet$] in ${\mathbb R}^{n}$ we fix the 1--norm: $ |I|:=|I|_1:=\sum_{1\leq i\leq n_1}|I_i|$;

\item[$\bullet$] in ${\mathbb T}^{n}$ we fix the ``sup--metric'': $  |\varphi|:=|\varphi|_{\infty}:=\max_{1\leq i\leq n}|\varphi_i|$ (mod $2\p$); 

\item[$\bullet$] in ${\mathbb R}$ we fix the sup norm: $  |(p, q)|:=|(p, q)|_{\infty}:=\max\{|p|, |q|\}$;

\item[$\bullet$] for  matrices  we use the ``sup--norm'': $    |\b|:=|\b|_{\infty}:=\max_{i,j}|\b_{ij}|$;

\item[$\bullet$] we denote as $B^n_{\varepsilon}(z_0)$ the complex ball having radius $\varepsilon$ centered at $z_0\in {\mathbb C}^n$. If $z_0=0$, we simply write $B^n_{\varepsilon}$.

\item[$\bullet$] if $A\subset {\mathbb R}^{n}$, and $r>0$, 
we denote by $A_r:=\bigcup_{x_0\in A} B^n_r(x_0)$
the complex $r$--neighbrhood of $A$ (according to the prefixed  norms/metrics above);

\item[$\bullet$] given $A\subset {\mathbb R}^n$ and positive numbers $r$, 
$\varepsilon$, 
$s$, we let
\beqano v:=(r, \varepsilon)\ ,\quad U_{v}:={A}_r
\times B^2_{ 
\varepsilon}\ ,\quad W_{v, s}:={U}_v
\times{\mathbb T}^{n}_{
s}\eeqano
\item[$\bullet$]  if $f$ is real--analytic on a complex domain of the form 
$W_{v_0, s_0}$, with $v_0=(r_0, \varepsilon_0)$, $r_0>r$, $\varepsilon_0>\varepsilon$, $s_0>s$, we denote by $\|f\|_{v,s}$ its
 ``sup--Taylor--Fourier norm'':
\beq{norm}\|f\|_{v, s}:=\sum_{k,\a,\b
 }\sup_{U_v}|f_{\a,\b, k}|e^{|k|s}\varepsilon^{|(\a,\b)|}\eeq with $|k|:=|k|_1$, $|(\a,\b)|:=|\a|_1+|\b|_1$,
where $f_{k, \a, \b}(I)$ denotes the  coefficients in the expansion 
\beqano f=\sum_{(k,\a,\b)\in \integer^{n}\times {\mathbb N}^\ell\times {\mathbb N}^\ell\atop{\a_i\ne \b_i\forall i}
  }f_{k, \a,\b}(I)e^{ik\cdot\varphi} p^\a q^\b\ ;\eeqano

\item[$\bullet$] if $f$ is as in the previous item, $K>0$  and ${\mathbb L}$ is a sub--lattice of $\integer^n$, $T_Kf$ and $\P_{\mathbb L} f$ denote, respectively, the $K$--truncation and the ${\mathbb L}$--projection of $f$:
$$T_Kf:=\sum_{(k,\a,\b)\in \integer^{n}\times {\mathbb N}^\ell\times {\mathbb N}^\ell\atop{\a_i\ne \b_i\forall i}
|k|_1\le K  }f_{k, \a,\b}(I)e^{ik\cdot\varphi}p^\a q^\b\ ,\quad \P_{\mathbb L} f:=\sum_{(k,\a,\b)\in \integer^{n}\times {\mathbb N}^\ell\times {\mathbb N}^\ell\atop{\a_i\ne \b_i\forall i}, k\in {\mathbb L}
  }f_{k, \a,\b}(I)e^{ik\cdot\varphi}p^\a q^\b$$
with $f_{k, \a,\b}(I):=f_{k, \a,\b}(I, 0, 0)$.
We say that $f$ {\it is $(K, {\mathbb L})$ in normal form} if $f=\P_{\mathbb L} T_{K}f$. If ${\mathbb L}$ is strictly larger than $\{0\}$, we say that $f$ {\it is resonant normal form}.

\end{itemize} 
\begin{proposition}[Partially hyperbolic averaging  theory]\label{average BCV}
Let 
 $H=h(I_1, I_2, pq)+        f(I,\f,p,q)
$ be a real--analytic function on $W_{v_0, s_0}$, with $v_0=(r_0, \varepsilon_0)$. Let $K$, 
$r$, $s$, $\varepsilon$, $\hat r$, $\hat s$,  positive numbers, with  
$\hat r<r/4$, 
 $\hat s<s/4$ and $\hat \varepsilon<\varepsilon/4$. 
 Put $\hat\sigma:=\min\left\{\hat s, \frac{\hat\varepsilon}{\varepsilon}\right\}$.
 Assume
 there exist positive numbers $\a_1$, $\a_2>0$, with $\alpha_1\ge \alpha_2$, such that,  for all $k=(k_1, k_2, k_3)\in{\mathbb Z}^{n+1}$, $0<|k|\le K$ and for all $(I, p, q)\in U_{r, \varepsilon}$, 
\beqa{hyper non res}
|\o_1\cdot k_1+\o_2\cdot k_2-\ii k_3\n |\ge \left\{
\begin{array}{lll}
\a_1\quad &{\rm if}\quad &k_1\ne 0\\
\a_2 &{\rm if}& k_1=0,\ (k_2, k_3)\ne (0, 0)\\
\end{array}
\right.
\eeqa
and
            \beqa{new smallness cond} 
        K\widehat \s\ge 8\log 2\,,\qquad   \frac{2^{ 3}c_1
           K \widehat \s}{\a_2 \d} \|f\|_{r, s, \varepsilon}<1,\qquad \d:= \min\{ \hat r \hat s,\ \hat\varepsilon^2
           \}  \eeqa
    with a suitable number $c_1$.     Then, one can find a real-analytic and symplectic transformation \[ \Phi_*:\quad W_{r_*, s_*, \varepsilon_*
    }\to W_{r, s, \varepsilon} \] 
    with $r_*=r-4\hat r$, $s_*=s-4\hat s$, $\varepsilon_*=\varepsilon-4\hat \varepsilon$,
     which conjugates $H$ to
    \begin{eqnarray*}
        H_{*}(I, \f,p,q):=H\circ \Phi_*={\rm h}(I,pq)+
        g(I,\f,p,q)+
        f_{*}(I,\f,p,q),
    \end{eqnarray*}
    where $g$ is $(K, \{0\})$ in normal form, and $g$, $f$ verify
     \beqa{results} &&\|g-\P_{0}T_Kf\|_{r_*, s_*, \varepsilon_*}\leq \frac{8c_1\,\|f\|^2_{r, s, \varepsilon}}{\a_2 \d}\nonumber\\
 &&     \|f_*\|_{r_*, s_*, \varepsilon_*}\le e^{-K\hat \s/4}\| f\|_{r_*, s_*, \varepsilon_*} \eeqa Finally, $\Phi_*$ verifies
 \beq{results2}
\max\big\{ \a_1\hat s|I_1-I_1'|\,,\, \a_2\hat s|{I_2}-{I_2}'|\,, \,  \a_2 \hat r\,|\varphi-\varphi'|\,, \, \a_2\hat \varepsilon\,|p-p'|\,, \, \a _2 \hat\varepsilon\,|q-q'| \big\}\,\leq2 { c_1E}\,.
\eeq
\end{proposition}
Proposition~\ref{average BCV} is an extension of the Normal Form Lemma by J. P\"oschel~\cite{poschel93}. The extension pertains at introducing the $(p, q)$ coordinates in the integrable part and leaving the amounts of analyticity $\hat r$, $\hat s$ and $\hat\varepsilon$ as independent. This is needed in order to construct the motions~\equ{whiskered tori}, where the coordinates $(\p, \k)$ are not set to $(0, 0)$, but take value in a small neighborhood of it. 
A more complete statement implying Proposition~\ref{average BCV} is quoted and proved in
 Section~\ref{sec: averageELLIPTIC}.

  \vskip.1in
  \noindent
  Below, we let $B:=B^2_{\ovl\varepsilon}(0)$; therefore, $B_\varepsilon$ will stand for $B^2_{\ovl\varepsilon+\varepsilon}(0)$. 
  
  \begin{lemma}[{\sc kam} Step Lemma]\label{step KAM}
Under the same assumptions and notations as in Theorem~\ref{two scales KAM},  there exists a sequence of numbers $\r_j$, $\varepsilon_j$, $s_j$; of domains $$(W_j)_{\rho_j, \varepsilon_j, s_j}=(A_j)_{\rho_j, \varepsilon_j}\times
{\mathbb T}^n_{\overline s+s_j}
\,,\qquad {\rm with}\quad (A_j)_{\rho_j, \varepsilon_j}:=\bigcup_{(p_j, q_j)\in B_{\varepsilon_j}}\left(D_j(p_jq_j)\right)_{\rho_j}\times\{(p_j, q_j)\}$$  and a real--analytic and symplectic transformations \beqa{Psi+}\Psi_{j+1}:\ (I_{j+1},\varphi_{j+1}, p_{j+1}, q_{j+1})\in (W_{j+1})_{\rho_{j+1}, \varepsilon_{j+1}, s_{j+1}} \ \to\ (I_{j},\varphi_{j}, p_{j}, q_{j})\in(W_j)_{\rho_j, \varepsilon_j, s_j}\eeqa such that 
\beqano
\HH_{j+1}(I_{j+1}, \varphi_{j+1}, p_{j+1}, q_{j+1})&=&\HH_j\circ \Psi_{j+1}(I_{j+1}, \varphi_{j+1}, p_{j+1}, q_{j+1})\nonumber\\
&=&{\rm h}_{j+1}(I_{j+1},  p_{j+1}q_{j+1})+{\rm f}_{j+1}(I_{j+1}, \varphi_{j+1}, p_{j+1}, q_{j+1})\eeqano and such that the following holds. Letting $E_0:=E$, $(M_{0}, \ovl M_{0}, \widehat M_{0}, L_{0})=(M, \ovl M, \widehat M, L)$, $s_0:=s$, $\rho_0:=\rho$, $\varepsilon_0:=\varepsilon_0$, $\lambda_0:=\lambda$ and, given,
for $0\le j\in {\mathbb Z}$, $E_j$, $(M_{j}, \ovl M_{j}, \widehat M_{j},  L_{j})$, $s_j$, $\rho_j$, $\varepsilon_j$, $\lambda_j$, define
%
\beqa{K} K_{j}&:=&\frac{32}{s_{j}}\log_+\Big(\frac{E _{j}L_{j}M_{j}^2}{\fg^2}\Big)^{-1}\\
\label{hat r}\widehat \r_{j}&:=&
\min\left\{
\frac{\g_1}{2M_{j}K_{j}^{\t+1}}\,,\ 
\frac{\g_2}{2\widehat  M_{j}K_{j}^{\t+1}}\,,\
\frac{\lambda_j}{2  M_{j}K_{j}}\,,\
\frac{\lambda_j}{2  \widehat M_{j}K_{j}}\,,\ 
\r_{j}\right\}\ ,
\\
\label{tilderhoj}
\widetilde\rho_j&:=&\min\left\{\widehat\rho_j, \frac{\varepsilon^2_j}{s_j}\right\}\,,\quad\widehat  E_{j}:=\frac{E_{j}L_{j}}{\widehat \r_{j}\widetilde \r_{j}}
\\
E_{j+1}&:=&\frac{E _{j}L_{j}M_{j}^2}{\fg^2}\ ,\quad (M_{j+1}, \ovl M_{j+1}, \widehat M_{j+1}, L_{j+1})=2(M_{j}, \ovl M_{j}, \widehat M_{j},  L_{j})\nonumber\\
\rho_{j+1}&:=&\frac{\widehat\rho_j}{4}\,,\ \varepsilon_{j+1}:=\frac{\varepsilon_j}{4}\,,\ \,,\lambda_{j+1}:= \lambda_j- 2^{8}\frac{E_j}{\varepsilon^2_j}\,,\quad s_{j+1}:=\frac{s_j}{4}\,.\nonumber\eeqa
Then, for all $(p_{j+1}, q_{j+1})\in B_{\varepsilon_{j+1}}$,
\begin{itemize}
\item[(i)] $D_{j+1}(p_{j+1}q_{j+1})\subseteq {(D_{j}}(p_{j+1}q_{j+1}))_{\widehat \r_{j}/ {4}}$. Letting $$\varpi_{j+1}:=\partial_{(I_{j+1}, p_{j+1}q_{j+1})}{\rm h}_{j+1}(I_{j+1}, p_{j+1}q_{j+1}))=(\omega_{j+1}(I_{j+1}, p_{j+1}q_{j+1}), \nu_{j+1}(I_{j+1}, p_{j+1}q_{j+1}))$$ the  map  {$I_{j+1}\to \omega_{j+1}(I_{j+1}, p_{j+1}q_{j+1})$ } is a diffeomorphism of $\big(D_{j+1}(p_{j+1}q_{j+1})\big)_{\r_j}$ verifying  {$$\omega_{j+1}(D_{j+1}(p_{j+1}q_{j+1})), p_{j+1}q_{j+1})=\omega_{j}(D_{j}(p_{j+1}q_{j+1})), p_{j+1}q_{j+1})\,.$$} The map 
\beqano
\widehat \iota_{j{+1}}(p_{j+1}q_{j+1})&=&(\widehat \iota_{j{+1,} 1}(p_{j+1}q_{j+1}),\widehat \iota_{j{+1,}2}(p_{j+1}q_{j+1})):\nonumber\\
D_{j}(p_{j+1}q_{j+1})&\to& D_{j+1}(p_{j+1}q_{j+1})\nonumber\\
 I_j(p_{j+1}q_{j+1})
&\to& I_{j+1}(p_{j+1}q_{j+1}):=\omega_{j+1}^{-1}\big(\omega_{j}(I_j, p_{j+1}q_{j+1}), p_{j+1}q_{j+1}\big)
\eeqano
verifies
\begin{eqnarray}\label{estimates for lj}
&&\sup_{{ D}_{j}}|\widehat \iota_{j+1, 1}(p_{j+1}q_{j+1})-\id|\leq3{n} \frac{\overline M_1}{\overline M}\widehat  E _{j}\widetilde \r_{j}\leq3{n} \widehat  E _{j}\widetilde \r_{j}\ ,\nonumber\\
&&\sup_{{ D}_{j}}|\widehat \iota_{j+1, 2}(p_{j+1}q_{j+1})-\id|\leq3{n} \frac{\overline M_2}{\overline M}\widehat  E _{j}\widetilde \r_{j}\leq3{n}\widehat  E _{j}\widetilde \r_{j}\\
\label{estimates for lj2}&&{\cal L}(\widehat \iota_{j+1}(p_{j+1}q_{j+1})-\id)\leq2^9{n} \widehat  E _{j}
\end{eqnarray}
\item[(ii)] the perturbation $\pert_j$ has sup--Fourier norm
\beqano\|\pert_j\|_{(W_j)_{\rho_j, \varepsilon_j, s_j}}\leq E_j
\eeqano
\item[(iii)] the real--analytic symplectomorphisms $\Psi_{j+1}$ in~\equ{Psi+} verify 
\beqa{bounds22}
&&\sup_{(W_{j+1})_{\r_{j+1}, \varepsilon_{j+1},s_{j+1}}}|{I_{j, 1}}(I_{j+1},\varphi_{j+1}, p_{j+1}, q_{j+1})-{I_{j+1,1}}|\leq\frac{3}{ {4
}} \frac{\widehat M_j }{ {M_j}}  \widehat  E_j{\widetilde \r_j}\nonumber\\
&&\sup_{(W_{j+1})_{\r_{j+1}, \varepsilon_{j+1},s_{j+1}}}|{I_{j, 2}}(I_{j+1},\varphi_{j+1}, p_{j+1}, q_{j+1})-{I_{j+1,2}}|\leq\frac{3}{ {4
}}\widehat E_j {\widetilde \r_j}\nonumber\\
&& \sup_{(W_{j+1})_{\r_{j+1}, \varepsilon_{j+1},s_{j+1}}}|\varphi_{j}(I_{j+1},\varphi_{j+1}, p_{j+1}, q_{j+1})-\varphi_{j+1}|\leq\frac{3}{ {4
}}\widehat E_j   s_j\nonumber\\
&& {\sup_{(W_{j+1})_{\r_{j+1}, \varepsilon_{j+1},s_{j+1}}}|p_{j}(I_{j+1},\varphi_{j+1}, p_{j+1}, q_{j+1})-p_{j+1}|\leq\frac{3}{ {4
}}\widehat E_j  \varepsilon_j}\nonumber\\
&& {\sup_{(W_{j+1})_{\r_{j+1}, \varepsilon_{j+1},s_{j+1}}}|q_{j}(I_{j+1},\varphi_{j+1}, p_{j+1}, q_{j+1})-q_{j+1}|\leq\frac{3}{ {4
}}\widehat E_j  \varepsilon_j}\,.
\eeqa
The rescaled dimensionless map $\check\Phi_{j+1}:=\id+{\mathbb 1}_{\widehat \r_0^{-1},s_0^{-1}, \varepsilon_0^{-1}}\left(\Phi_{j+1}-\id\right)\circ {\mathbb 1}_{\widehat \r_{0},s_{0}, \varepsilon_{0}}$
has Lipschitz constant
  on $(W_{j+1})_{\r_{j+1}/\widehat\rho_0, \varepsilon_{j+1}/\varepsilon_0,s_{j+1}/s_0}$ 
\begin{eqnarray}\label{Psij close to id}
{\cal L}(\check\Phi_{j+1}-\id)\leq  6(n+1)\big(12\cdot(24)^\tau\big)^{j}\widehat  E _{j}\,;\end{eqnarray}
\item[(iv)] for any $j\ge 0$,  $\widehat  E _{j+1}<\widehat  E _{j}^2$, $\lambda_j\ge \frac{\lambda_0}{2}$.\end{itemize}
\end{lemma}

\vskip.1in
\noindent 
{\bf Proof}\ The proof of this proposition is obtained  generalizing~\cite[Lemma B.1]{chierchiaPi10}. We shall limit ourselves to describe only the different points, leaving to the interested reader the easy work of completing details.

\noindent
We construct the transformations~\equ{Psi+}
 by recursion,   based on  Proposition~\ref{average BCV}. 
  For simplicity of notations, we  shall 
 systematically eliminate the sub--fix ``$j$'' and replace ``$j+1$'' with a ``+''. As an example,
 instead of~\equ{Psi+},
 we shall write
 $$\Psi_{+}:\ W_{+}\to W\,.$$
 When needed, the base step will be labeled as ``$0$'' (e.g,~\equ{step0} below).
Let us assume (inductively) that  
 \beqa{good freq}&& \omega(D, pq)\subset{\cal D}_{\gamma_1, \gamma_2, \tau}\quad \forall \ (p, q)\in B_\varepsilon
 \\
 \nonumber\\
 \label{KAM ass}&&\widehat  c\widehat  E<1\\\nonumber\\
 \label{lambdanotsmall}&& \lambda\ge \max\left\{\frac{\gamma_2}{K^{\tau}}\,,\ \frac{\lambda_0}{2}\right\}\,.
 \eeqa
 Condition~\equ{good freq} is verified at the base step provided one takes $D_0=\omega_0^{-1}({\cal D}_{\gamma_1, \gamma_2, \tau}, p_0q_0)$; 
\equ{KAM ass} is so by assumption, while~\equ{lambdanotsmall} follows from~\equ{lambda cond}:
\beqa{lambdastep0}\lambda_0\ge \frac{\lambda_0}{2}\ge \frac{s_0^\tau\gamma_2}{6^{\tau}} \ge  \frac{\gamma_2}{K_0^{\tau}}\,.\eeqa
We aim to apply Proposition~\ref{average BCV} with  $\varepsilon$, $s$ of Proposition~\ref{average BCV} corresponding now 
to  $\ovl\varepsilon+\varepsilon$, $\ovl s+s$, and
 $$r=\widehat \rho\,,\quad \hat r=\frac{\hat\rho}{8}\,,\quad \hat s:=\frac{s}{8}\,,\quad \hat \varepsilon:=\frac{\varepsilon}{8}\,,\quad {\mathbb L}=\{0\}\,.$$
We  check that that~\equ{good freq} and~\equ{KAM ass} imply conditions~\equ{hyper non res} and~\equ{new smallness cond}. We start with~\equ{hyper non res}. If $(I, p, q)=(I_{1}, I_{2}, p_, q)\in A_{\widehat\rho, \varepsilon}$ and $k\in \integer^3\setminus\{0\}$, with $|k|_1\le K$, then there exists some $I_0(pq)=(I_{01}(pq), I_{02}(pq))$ such that $|I-I_0(pq)|<\widehat\rho$ and $\omega(I_0(pq), pq)=(\omega_{01}, \omega_{02})\in {\cal D}_{\gamma_1, \gamma_2, \tau}$. We have
\beqano
|\varpi(I,  pq)\cdot k|&=&\Big|\o_{01}\cdot k_1+\o_{02}\cdot k_2+(\o_1(I,  pq)-\o_1(I(pq),  pq))\cdot k_1\nonumber\\
&&+(\o_2(I,  pq)-\o_2(I(pq),  pq)))\cdot k_2-\ii\n(I,  pq) k_3 \Big|\nonumber\\
\nonumber\\
&\ge& \left\{
\begin{array}{lll}
\displaystyle\min\left\{\frac{\gamma_1}{2K^\tau}\,,\ \frac{\lambda}{2}\right\}\quad &{\rm if}\quad &k_1\ne 0\\\\
\displaystyle\min\left\{\frac{\gamma_2}{2K^\tau}\,,\ \frac{\lambda}{2}\right\} &{\rm if}& k_1=0,\ k_2\ne 0\\\\
\displaystyle\lambda&{\rm if}& k_1=k_2=0,\ k_3\ne 0
\end{array}
\right.\nonumber\\
&\ge& \left\{
\begin{array}{lll}
\displaystyle\alpha_1:=\frac{\gamma_1}{2K^\tau}\quad &{\rm if}\quad &k_1\ne 0\\\\
\displaystyle\alpha_2:=\frac{\gamma_2}{2K^\tau} &{\rm if}& k_1=0,\ (k_2, k_3)\ne (0, 0)
\end{array}
\right.
\eeqano
having used~\equ{lambdanotsmall}.
The bounds above have been obtained considering separately the cases $k_3\ne0$ and $k_3= 0$, and:

\noindent
-- if $k_3\ne 0$, taking the infimum of the modulus of the imaginary part of the expression between the $|$'s;
observing that $\ovl\omega_0=(\ovl\omega_{01}, \ovl\omega_{02})$ are real and bounding the differences $|\Im\big(\o_i(I,  pq)-\o_i(I(pq),  pq)\big)|$ with $MK\widehat\rho$ (when $i=1$), $\widehat MK\widehat\rho$ (when $i=2$) and using the definition of $\widehat\rho$ in~\equ{hat r}.

\noindent
--  if $k_3=0$, using the Diophantine inequality and again  bounding   the differences $|\Im\big(\o_i(I,  pq)-\o_i(I(pq),  pq)\big)|$ as in the previous case and using  the definition $\widehat\rho$.

\noindent
We now check condition~\equ{new smallness cond}. 
The inequality $Ks>  8\log 2$ 
is trivial by definition of $K$ (see~\equ{K}) and also the smallness condition ~\equ{new smallness cond}  is easily met, since $\hat\sigma=\min\{\frac{1}{8}\frac{\varepsilon}{\ovl\varepsilon+\varepsilon}\,,\ 
\frac{s}{8}
\}=\frac{s}{8}$, $\delta=2^{-6}\min\{\hat\rho s\,,\ \varepsilon^2\}=2^{-6}\widetilde \r s$ (by the definition of $\widetilde\rho$ in~\equ{tilderhoj}), whence
$$2^{3} c_1\frac{K\frac{s}{8}}{\a_2 \d}\|f\|_{ {W_{\widehat \r, \varepsilon, s}}}\leq 2^{6}c_1\frac{E L}{\widehat \r\widetilde \r}\le \widehat  c \widehat  E<1$$
having used $L\ge \widehat  M^{-1}$,  $M^{-1}$, so $\a_2\geq K L^{-1}\widehat \r$,  $2^{6}c_1<\widehat  c$,  and~\equ{KAM ass}. Thus, by Proposition~\ref{average BCV}, $\HH$ may be conjugated to
$$\HH_+:=\HH\circ \Psi_+=\hh_+(I_+, p_+q_+)+f_+(I_+, \f_+, p_+, q_+)$$
where
$$\nfj+(I_+, p_+q_+)=\nf(I_+, p_+q_+)+{\rm g}(I_+, p_+q_+)$$
while,
by ~(\ref{results}) and the choice of $K$,
\beqa{F++}
\|f_+\|_{{\widehat \r/2, \varepsilon/2, s/2}}&\leq& e^{-Ks/32}E\leq \frac{E LM  ^2}{\fg^2}E=E_+\ .\eeqa
The conjugation is realized by  an analytic transformation
\beqno
\Psi_+:\quad(I_+,\varphi_+, p_+, q_+,)\in W_{\widehat \r/2, \varepsilon/2, s/2}\to (I,\varphi, p, q)\in W_{\widehat \r, \varepsilon, s}\ .
\eeqno
Using~(\ref{results2}), $\widetilde\rho\le \varepsilon^2/s$,  $\a_1\geq  {MK\widehat \r}$,  $\a_2=\frac{\gamma_2}{2K^\tau}\geq  {L^{-1} K\widehat \r}$, $Ks\geq  {6}$  {and the definition of $\widehat E$}, we obtain the bound~\equ{bounds22} with, at the left hand side, the set $W_{\widehat \r/2, \varepsilon/2, s/2}$. Below we shall prove that 
${W_{+}}_{\r_{+}, \varepsilon_{+},s_{+}}\subset W_{\widehat \r/2, \varepsilon/2, s/2}$, so we shall have~\equ{bounds22}.\\
We now evaluate the generalized frequency
$$\varpi_+(I_+, p_+q_+):=\partial_{I_+, p_+q_+}{\rm h}_+(I_+, p_+q_+)=\big(\omega_+(I_+, p_+q_+)\,,\ \nu_+(I_+, p_+q_+)\big)\,.$$
with
\beqa{omega+}\omega_+(I_+, p_+q_+):=\partial_{I_+}{\rm h}_+(I_+, p_+q_+)=\partial_{I_+}{\rm h}(I_+, p_+q_+)+\partial_{I_+}{\rm g}(I_+, p_+q_+)\eeqa
(the ``new frequency map'') and
\beqa{nu+} \nu_+(I_+, p_+q_+)&:=&\partial_{p_+q_+}{\rm h}_+(I_+, p_+q_+)=\nu(I_+, p_+q_+)+\partial_{p_+q_+}{\rm g}(I_+, p_+q_+)\eeqa
(the ``new Lyapunov exponent''). 
\begin{lemma}\label{new integrable part}
Let $(p_+, q_+)\in B_{\varepsilon/2}$. The new frequency map $ {\o}_+$ is injective on $ D(p_+q_+)_{\widehat \r/2}$ and maps  $ D(p_+q_+)_{\widehat \r/4}$ over $ {\o}(D, p-+q_+)$. The map  $\widehat \iota_+(p_{+}q_{+})=(\widehat \iota_{+1}(p_{+}q_{+}),\widehat \iota_{+2}(p_{+}q_{+})):= {\o}_+^{-1}\circ {\o}|_{D(p_+q_+)}$ which assigns to a point $I_0\in D(p_+q_+)$ the $ {\o}_+(\cdot, p_+q_+)$--preimage  of $ {\o}(I_0, p_+q_+)$  in ${ D}(p_+q_+)_{\widehat \r/4}$ satisfies
\begin{eqnarray}\label{lip}
&& \sup_{(A_+)_{\rho_+, \varepsilon_+}}|\widehat \iota_{+1}(p_{+}q_{+})-\id|\leq {3}{n}\frac{\inHs_1  E}{\widehat \r}\leq {3}{n}\frac{\inHs  E}{\widehat \r}\ ,\nonumber\\
&&\sup_{(A_+)_{\rho_+, \varepsilon_+}}|\widehat \iota_{+2}(p_{+}q_{+})-\id|\leq3{n}\frac{\inHs_2  E}{\widehat \r}\leq3{n}\frac{\inHs  E}{\widehat \r}\ ,\nonumber\\
&& {\cal L}(\widehat \iota_+(p_{+}q_{+})-\id)\leq2^{ {9}}{n}\frac{\inHs E}{\widehat \r^2}\ .
\end{eqnarray}
The Jacobian matrix $U_+:=\partial^2_{ {I_+}}\nfj+ {(I_+, p_+q_+)}$ is non singular on $ D_{\widehat \r/ {4}} {\times B^2_{\varepsilon/2}}$ and the following bounds hold 
\beqa{bounds on M}
&&M_+ :=2M\ge \sup_{(A_+)_{\rho_+, \varepsilon_+}}\|U_+\|\ ,\quad  \widehat M_+:=2\widehat M\ge \sup_{(A_+)_{\rho_+, \varepsilon_+}}\|\widehat  U_+\|\ ,\nonumber\\
&& \ovl M_+:=2\ovl M\ge \sup_{(A_+)_{\rho_+, \varepsilon_+}}\|U_+^{-1}\|\ ,\quad\ovl M_{i+}:=2\ovl M_i\ge \sup_{(A_+)_{\rho_+, \varepsilon_+}}\|T_{i+}\|\ ,\quad i=1,\ 2\ .\eeqa
where
$ U_+^{-1}=:\left(\begin{array}{lrr} T_{+1}\\ T_{+2}
\end{array}
\right)$.
Finally, the new Lyapunov exponent $\nu_+(I_+, p_+q_+)$
satisfies \beqa{lambda+}\l_+:=\l- 2^4\frac{E}{\varepsilon^2}\le  \inf_{(A_+)_{\rho_+, \varepsilon_+}}|\Re\n_+|\,.\eeqa
\end{lemma}
Postponing for the moment the proof of this Lemma, we let $\r_+:=\widehat \r/ 2$, $s_+:=s/2$,  {$\varepsilon_+=\varepsilon/2$} and $ D_+(p_{+}q_{+}):=\widehat \iota_+(p_{+}q_{+})(D(p_{+}q_{+}))$. By Lemma~\ref{new integrable part}, $D_+$ is a subset of $D_{\widehat \r/4}$ and hence
\beq{I+r+}(D_+)_{\r_+}\subset  D_{\widehat \r/2}\ .\eeq 
We  prove  that  
$
\widehat  E _+=\frac{E_+L_+}{\widehat \r_+^2}\leq \widehat  E ^2
$. Since
\beq{s+}
s_+=\frac{s}{4}\quad \textrm{and}\quad x_+:=\Big(\frac{E_+ L_+M _+^2}{\fg^2}\Big)^{-1}=\frac{x^2}{8}\quad \textrm{where}\quad x:=\Big(\frac{E LM ^2}{\fg^2}\Big)^{-1}
\eeq
we have
\beqa{K+}K_+=\frac{2^5}{s_+}\log x_+=
\frac{2^7}{s}\log \frac{x^2}{8}
=\frac{2^8}{s}\log_+ x-\frac{3\cdot 2^7}{s}\log_+ 2<8 K\,.\eeqa
Finally,~\equ{KAM cond},~\equ{lambda+}  and the definition of $\tilde E$ imply
$\lambda_+\ge \frac{\lambda}{2}$. Collecting all bounds, we get
\beqa{r+}
\widehat \r_+&=&
\min\left\{
\frac{\g_1}{ 2M_{+}K_{+}^{\t+1}}\,,,\ 
\frac{\g_2}{ 2\widehat  M_{+}K_{+}^{\t+1}}\,,\
\frac{\lambda_+}{ 2  M_{+}K_{+}}\,,\
\frac{\lambda_+}{ 2  \widehat M_{+}K_{+}}\,,\ 
\r_{+}=\frac{\widehat\rho}{2}\right\}\geq\frac{\widehat \r}{2\cdot 8^{\t+1}}\nonumber\\
\widetilde\rho_+&=&\min\left\{\widehat\rho_+, \frac{\varepsilon^2_+}{s_+}\right\}
\geq\frac{\widetilde \r}{2\cdot 8^{\t+1}}
\eeqa
and 
$$\widehat  E _+=\frac{E_+L_+}{\widehat \r_+\widetilde\rho_+}\leq \frac{E^2LM ^2}{\fg^2}\frac{2L}{\widehat \r\widetilde\rho}4\cdot 8^{2(\t+1)}=8\cdot 8^{2(\t+1)}\frac{E\,LM ^2}{\fg^2}\widehat  E $$
Now, using, in the last inequality, the bound 
$$\frac{E\,LM ^2}{\fg^2}
\le\frac{1}{4}\left(\frac{s}{6}\right)^{2(\t+1)}
\frac{E L}{\widehat \r^2}
 \le\frac{1}{4}\left(\frac{s}{6}\right)^{2(\t+1)}\widehat  E $$
(since $\hat\rho\le \frac{\g_1}{ 2M K^{\t+1}}$ and  $K\geq 6/s$) we find
\beq{KAM+}\widehat  E _+\leq 2(\frac{4}{3}s)^{\t+1}\widehat  E ^2<\widehat  E ^2\eeq
(having used $s\le1/2$).  We now prove that $\lambda_{+}\ge \frac{\lambda_0}{2}$. Iterating~\equ{lambda+} and using $\widehat\rho_k\le \widehat\rho_{k-1}/4$, $\widetilde\rho_k\le \widetilde\rho_{k-1}/4$, $\varepsilon_k= \varepsilon_{k-1}/4$, $L_k=2L_{k-1}$,~\equ{KAM+} and
the second condition in~\equ{KAM cond} with $\tilde c=2^6$,
 we get
\beqa{step0}
\lambda_{+}&=&\lambda_{j+1}=\lambda_0-2^4\sum_{k=1}^{j}\frac{E_k}{\varepsilon_k^2}\ge 
\lambda_0-2^4\sum_{k=1}^{j}\widehat E_k\frac{\widehat\rho_k\widetilde\rho_k}{\varepsilon^2_k L_k}
\ge 
\lambda_0-2^4\frac{\widehat\rho_0\widetilde\rho_0}{\varepsilon^2_0 L_0}\sum_{k=1}^{j}\widehat E_k\ge 
\lambda_0-2^5\frac{\widehat\rho_0\widetilde\rho_0}{\varepsilon^2_0 L_0}\widehat E_0\nonumber\\
&=&\lambda_0-2^5\frac{ E_0}{\varepsilon^2_0}
\ge \frac{\lambda_0}{2}\,.
\eeqa
This allows to check~\equ{lambdanotsmall} at the next step: using~\equ{lambdastep0} and~\equ{K+}, we have
$$\lambda_+\ge \frac{\lambda_0}{2}\ge   \frac{\gamma_2}{K_0^{\tau}} \ge  \frac{\gamma_2}{K_+^{\tau}}\,.$$
Finally,~\equ{estimates for lj} and~\equ{estimates for lj2} follow from~\equ{lip}, while
the estimate in~(\ref{Psij close to id}) is a consequence of ~(\ref{bounds22}),~(\ref{I+r+}),~(\ref{s+}),~(\ref{r+}), {inequality $LM\ge 1$} and Cauchy estimates:
\beqano
{\cal L}(\check{{\Phi}}_{j+1}-\id)&\leq&2(n+1)\sup_{{(\check W_{j+1})_{\r_{j+1}, \varepsilon_{j+1},s_{j+1}}}}\|D(\check{{\Phi}}_{j+1}-\id)\|_{\infty}\nonumber\\
&\leq&2(n+1) \frac{\frac{3}{4}\widehat  E_j\max\{\widehat\r_{j}/\r_0,s_j/s_0, {{\varepsilon_j}/{\varepsilon_0}}\}}{\min\{\widehat \r_{j}/(4\widehat \r_0),s_{j}/(4s_0), {{\varepsilon_{j}}/{(4\varepsilon_0)}}\}}\nonumber\\
&\leq& 2(n+1)\frac{3/4(1/4)^{j}}{1/4\left(\frac{1}{2({24})^{\t+1}}\right)^{j}}\widehat  E _{j}=6(n+1)\big(12\cdot(24)^\tau\big)^{j}\widehat  E _{j}\ .\qquad \square
\eeqano

\vskip.1in
\noindent
{\bf Proof of Lemma~\ref{new integrable part}}\ The proof of this proposition is obtained  generalizing~\cite[Lemma B.2]{chierchiaPi10}. As above, we limit to discuss only the different parts.

\noindent
By~(\ref{results}),
\begin{eqnarray*}
\sup_{ {{ D}_{\widehat \r/2}\times B^2_{\varepsilon/2}}}|{\rm g}|\leq \sup_{ {{ D}_{\widehat \r/2}\times B^2_{\varepsilon/2}}}|{\rm g}-\overline\pert|+\sup_{ {{ D}_{\widehat \r/2}\times B^2_{\varepsilon/2}}}|\overline\pert|\leq  {\frac{3}{2}}E\ ,
\end{eqnarray*}
(where $\overline f$ denotes the average of $f$). Therefore
we may bound
\begin{eqnarray*}
%
%
&&\sup_{ {{ D}_{\widehat \r/4}\times B^2_{\varepsilon/2}}}\|(\partial_{I_+}^2\,{\rm h})^{-1}\partial_{I_+}^2\,{\rm g}\|\leq 2 \ovl M \frac{\frac{3}{2}E}{(\widehat \r/4)^2}\leq 2^{ {6}}\frac{\ovl ME}{\widehat \r^2}\leq 2^{ {6}}\frac{\ovl ME}{\widehat \r^2}<\frac{1}{2}\ \nonumber\\
\end{eqnarray*}
This shows that the function~\equ{omega+}
has a Jacobian matrix
$$\partial_{I_+}\omega_+(I_+, p_+q_+)=\partial_{I_+}^2\,{\rm h}_+(I_+, p_+q_+)=\partial^2_{I_+}{\rm h}(I_+, p_+q_+)+\partial^2_{I_+}{\rm g}(I_+, p_+q_+)$$
which is invertible for all $(p_+, q_+)\in B^2_{\varepsilon/2}$ and satisfies
$$\ovl M_+:=\sup_{ {{ D}_{\widehat \r/4}\times B^2_{\varepsilon/2}}}\Big\|\Big(\partial_{I_+}\omega_+(I_+, p_+q_+)\Big)^{-1}\Big\|\le 2\ovl M$$
In a similar way one proves~\equ{bounds on M}.
Next, for any fixed $(p_+, q_+)\in B^2_{\varepsilon/2}$ and $\ovl\omega=\omega(I(p_+q_+), p_+q_+)\in \omega(D, p_+q_+)
$ with $I(p_+q_+)\in D$, we want to find $I_+=I_+(p_+q_+)\in D_+$ such that 
\beqa{implicit equation}\omega_+(I_+(p_+q_+), p_+q_+)=\ovl\omega=\omega(I(p_+q_+), p_+q_+)\eeqa
To this end, we
consider the function
\beqano I_+\in { D}_{\widehat \r/2}\to F(I_+, p_+q_+):=\omega_+(I_+, p_+q_+)-\ovl \omega \quad (p_+, q_+)\in B^2_{\varepsilon/2}\eeqano
As $F$ differs from $\omega_+$ by a constant, we have
$$m:=\sup_{ {{ D}_{\widehat \r/4}\times B^2_{\varepsilon/2}}}\Big\|\Big(\partial_{I_+}F(I_+, p_+q_+)\Big)^{-1}\Big\|=\sup_{ {{ D}_{\widehat \r/4}\times B^2_{\varepsilon/2}}}\Big\|\Big(\partial_{I_+}\omega_+(I_+, p_+q_+)\Big)^{-1}\Big\|\le 2M\,.$$
Similarly, we bound the quantities
$$Q:=|\partial^2_{I_+}F(I)|=|\partial^3_{I_+}{\rm g}(I_+, p_+q_+)|\le 6\frac{\frac{3}{2}E}{(\widehat\rho/4)^3}<2^{10}\frac{E}{\widehat\rho^3}.$$
and
$$P:=|F(I(p_+q_+))|=|\partial_{I_+}{\rm g}(I(p_+q_+), p_+q_+)|\le \frac{\frac{3}{2}E}{(\widehat\rho/4)}\le 2^{3}\frac{E}{\widehat\rho}\,.$$
Putting everything together, we get
$$4m^2 PQ\le 2^{16}\frac{M^2 E^2}{\widehat\rho^4}\le \widehat c^2 \widehat E^2<1$$
By the Implicit Function Theorem (e.g.,~\cite[Theorem 1 and Remark 1]{cellettiC98}), Equation~\equ{implicit equation} has a unique solution $$(p_+, q_+)\in B_{\varepsilon/2}\to I_+(p_+q_+)\in B_r(I(p_+q_+))\,,$$ with
$$r=2mP\le 2^5\frac{ME}{\widehat\rho}\le \frac{\widehat\rho}{4}$$
so we can take
$$D_+(p_+q_+)=\bigcup_{\ovl\omega\in\omega(D, p_+q_+)}\{I_+(p_+q_+)\}$$
This ensures that~\equ{good freq} holds also for $D_+$.\\
Finally, the real part of the function~\equ{nu+}
satisfies the lower bound
$$\inf_{ {{ D}_{\widehat \r/2}\times B^2_{\varepsilon/4}}}\left|\Re \n_+\right|\ge \l- \frac{E}{(\varepsilon/4)^2}=\l_+\,.$$
The proof of~\equ{lip} proceeds as in~\cite[proof of Lemma B.2]{chierchiaPi10}. $\qquad \square$
\vskip.1in
\noindent
{\bf Proof of Theorem~\ref{two scales KAM}}.

\noindent
{\bf Step 1} {\it Construction of the ``generelised limit actions''}

\noindent
Let $(\p, \k)\in B_0=B^2_{\ovl\varepsilon}=\bigcap_{j\ge 0}B_{\varepsilon_j}$. Define, on $D_0(\p\k)=\omega_0^{-1}(\cD_{\gamma_1, \gamma_2, \tau}, \p\k)\cap D$,
$$\check \iota_j(\p\k):=\widehat \iota_j(\p\k)\circ \widehat \iota_{j-1}(\p\k)\circ \cdots\circ  \widehat \iota_1(\p\k)\quad j\ge 1\,.$$
Then $\check\iota_j(\p\k)$ converge uniformly to a $\check\iota (\p, \k)=(\check\iota_1(\p, \k), \check\iota_2(\p, \k))$ verifying
\beqa{ell cl to id}\sup_{D_0(\p\k)}|\check\iota_1(\p\k)-\id|\le 6n \frac{\ovl M_1}{\ovl M}\widetilde\rho_0\widehat E_0\,,\quad \sup_{D_0(\p\k)}|\check\iota_2(\p\k)-\id|\le 6n \frac{\ovl M_i}{\ovl M}\widetilde\rho_0\widehat E_0\,.\eeqa
Moreover,
as
$$\sup |\widehat\iota_j(\p\k)-\widehat\iota(\p\k)|\le 6n \widehat E_j\widetilde\rho_j<\frac{6n}{\widehat c}\widehat\rho_j<\rho_j$$
we have
\beqa{intersection}D_*(pq):=\check \iota(\p\k) (D_0(\p\k))\subset \bigcap_{j} {D_j(\p\k)}_{\rho_j}\,.\eeqa
In particular, taking $j=0$,
\beqa{where is Cantor} D_*(\p\k)\subset (D_0(\p\k))_{6n \widehat E_0\widetilde\rho_0}\,.\eeqa
Moreover,
$${\cal L}(\check\iota(\p\k)-\id)\le 2^8 n \widehat E\,.$$
So $\check\iota(\p\k)$ is bi--Lipschitz, with
$${\cal L}_-(\check\iota(\p\k))\ge 1-2^8 n \widehat E\,,\qquad {\cal L}_+(\check\iota(\p\k))\le 1+2^8 n \widehat E\,.$$

\noindent
{\bf Step 2} {\it Construction of $\phi_{\omega_*}$.} 
For each $j\ge 1$, the transformation $$\Phi_j:=\Psi_1\circ\cdots\circ \Psi_j$$
is defined on $(W_j)_{\rho_j, s_j, \varepsilon_j}$. If
$$A_*:=\bigcup_{|(\p, \k)|<\ovl\varepsilon}D_*(\p\k)\times\{(\p, \k)\}\,,\qquad W_*:=A_*\times {\mathbb T}^n\,.$$
then, by~\equ{intersection}, $W_*\subset \bigcap_j (W_j)_{\rho_j, s_j, \varepsilon_j}$. The sequence $\Phi_j$ converges uniformly on $W_*$ to a map $\Phi$. We then let
\beqano
\phi_{\omega_*}(\vartheta, \p, \k)&=&\Big(v(\vartheta, \p, \k; \omega_*), \vartheta+u(\vartheta, \p, \k; \omega_*), 
\p+w(\vartheta, \p, \k; \omega_*), \k+y(\vartheta, \p, \k; \omega_*)
\Big)\nonumber\\
&:=&\Phi\left(\check\iota(\omega_0^{-1}(\omega_*, \p\k)), \vartheta, \p, \k \right)
\eeqano
with $v(\vartheta, \p, \k; \omega_*):=\big(v_1(\vartheta, \p, \k; \omega_*), v_2(\vartheta, \p, \k; \omega_*)\big)$. 
  Since~(\ref{bounds22}) imply, on $W_*$\footnote{$\P_z$ denotes the projection on the $z$--variables.},
  \beq{bounds on Phi11}\sup_{W_*}|\P_{I_1}\Phi-\id|_1\leq2n\frac{\sHs_0}{\Hs_0}\widehat E_0\widetilde\r_0\eeq
and similarly,
\beqa{bounds on Phi2}
&&\sup_{W_*}|\P_{I_2}\Phi-\id|_1\leq 2n \widehat E_0\widetilde\r_0\ ,\quad \sup_{W_*}|\P_{\varphi}\Phi-\id|_{\infty}\leq 2\widehat E_0 s_0\ ,\nonumber\\
&& \sup_{W_*}|\P_{p}\Phi-\id|_{\infty}\leq 2\widehat E_0 \varepsilon_0\,,\quad  \sup_{W_*}|\P_{q}\Phi-\id|_{\infty}\leq 2\widehat E_0 \varepsilon_0
\eeqa
then, in view of~(\ref{ell cl to id}),~(\ref{bounds on Phi11}),~(\ref{bounds on Phi2}), the definition of $W_*$ and the triangular inequality, we have~\equ{***}.
Equations~(\ref{where is Cantor}),~(\ref{bounds on Phi11}),~(\ref{bounds on Phi2}) also imply
\beqa{To*}{\rm T}_{\o_*}:=\phi_{\o_*}(\torus^n, 0, 0)\subset  (D_*(0))_{2 \widehat E_0\widetilde\r_0}\times\torus^n\times B^2_{r'}\subset (D_0(0))_r\times\torus^n\times B^2_{r'}\eeqa
where
$$r=8 n\widehat E_0\widetilde\r_0\,,\qquad r'=2\widehat E_0\varepsilon_0$$
Finally, with similar arguments as in {Step 1}, by~(\ref{LipcheckPhi}), the rescaled map  $$\check\Phi:=\id+{\mathbb 1}_{\widehat \r_0^{-1},s_0^{-1}, \varepsilon_0^{-1}}\left(\Phi-\id\right)\circ {\mathbb 1}_{\widehat \r_{0},s_{0}, \varepsilon_{0}}$$ has Lipschitz constant
\beq{LipcheckPhi}{\cal L}(\check\Phi-\id)\leq 2^6(n+1)  \widehat E_0\ .\eeq 
In particular, $\check\Phi$, hence, $\Phi$, and, finally, the map $(\vartheta,\p, \k; \o)\to \phi_\o(\vartheta, \p, \k)$ are bi--Lipschitz, hence, injective.
\vskip.1in
\noindent
{\bf Step 3} {\it For any $\o_*\in\cD_{\gamma_1, \gamma_2,\t}\cap\o_0(D, 0)$,  ${\rm T}_{\o_*}$ in~\equ{To*} is a $n$--dimensional $\rm H$--invariant torus with frequency $\o_*$.}
This assertion  is a trivial generalization of its analogue one in ~\cite[Proof of Proposition 3, Step 3]{chierchiaPi10}, therefore its proof is omitted.

\vskip.1in
\noindent
{\bf Step 4} {\it Measure Estimates (proof of~(\ref{tori measure}))} The proof of~(\ref{tori measure}) proceeds as in~\cite[Proof of Proposition 3, Step 4]{chierchiaPi10}, just replacing the quantities that in~\cite[Proof of Proposition 3, Step 4]{chierchiaPi10} are called
$$D_0\,,\quad D_*\,,\quad \check\iota\,,\quad \check\Phi\,,\quad {\rm K}$$
with the quantities here denoted as
$$D_0(0)\,,\quad D_*(0)\,,\quad \check\iota(0)\,,\quad \check\Phi\big|_{(\p, \k)=(0, 0)}\,,\quad {\rm K}_0\,.\qquad \square$$

\subsection{Normal Form Theory}\label{sec: averageELLIPTIC}
Proposition~\ref{average BCV} can be obtained from the more  general Proposition~\ref{averageELLIPTIC} below,  taking $m=1$, ${\mathbb L}=\{0\}$ and  changing coordinates as follows
\beqano p=\frac{{p_1}-\ii {q_1}}{\sqrt2}\ ,\qquad q=\frac{{p_1}+\ii {q_1}}{\sqrt2 \ii}\ .\eeqano
We define $c_m$ to be the smallest number such that, for any two functions, real--analytic in $W_{r, s, \varepsilon}$ and any choice of $\hat r<r$, $\hat s<s$, $\hat\varepsilon<\varepsilon$, 
$$\|\{f, g\}\|_{r-\hat r, s-\hat s, \varepsilon-\hat\varepsilon}\le \frac{c_m}{\delta}\|f\|_{r, s, \varepsilon}\|g\|_{r, s, \varepsilon}\quad {\rm with}\ \delta:=\min\{\hat r\hat s, \hat\varepsilon^2\}\,.$$

\begin{proposition}\label{averageELLIPTIC}
Let $\{0\}\subset {\mathbb L}\subset {\mathbb Z}$.
Proposition~\ref{average BCV} holds true taking
 $$H(I,\f,p,q)=h\left(I_1, I_2, J(p, q)\right)+ f(I,\f,p,q)\,,\quad       J(p, q):=\left( \frac{p_1^2+q_1^2}{2}\,,\ldots, \frac{p_m^2+q_m^2}{2}\right)
$$ 
replacing  $c_1$ with $c_m$, $\P_0$ with $\P_{\mathbb L}$ and  condition~\equ{hyper non res} with

\beqa{non res}
&&|\o_1\cdot k_1+\o_2\cdot k_2 |\ge \left\{
\begin{array}{lll}
\a_1\quad &{\rm if}\quad &k_1\ne 0\\
\a_2 &{\rm if}& k_1=0,\ k_2\ne 0\\
\end{array}
\right.\\\nonumber\\
&& \forall\ k=(k_1, k_2)\in{\mathbb Z}^{n_1}\times {\mathbb Z}^{n_2+m}\setminus{\mathbb L}\ne (0, 0)\,,\ |k|_1\le  K\,,\quad \forall (I_1, I_2, p, q)\in V_r\times B^{2m}_{\varepsilon}\nonumber
\eeqa
where
\beqano
\omega&=&(\omega_1, \omega_2)\nonumber\\
&:=&\left(\partial_{I_1} h\left(I_1, I_2, J(p, q)\right)\,,\ \partial_{\left(I_2,  J(p, q)\right)} h\left(I_1, I_2,  J(p, q)\right)
\right)\,.
\eeqano
\end{proposition}

 \begin{lemma}
    \label{itera1} Let $\hat r<r/2$ $\hat s<s/2$, $\hat\varepsilon<\varepsilon/2$ and $\d:=\min\{\hat r\hat s,\ \hat\varepsilon^2\}$. Let \[ H(u,\varphi,p,q)={\rm h}(I, p,q)+g(u, \f, p,q) +f(u, \f,p,q)\qquad g(u, \f, p,q)=\sum_{i=1}^mg_i(u, \f, p,q)\] be real-analytic on $W_{v,s,\varepsilon}$. Assume that inequality~\equ{non res} and \beqano \|f\|_{v,s, \varepsilon}<\frac{\alpha_2 \d}{c_m} \eeqano are satisfied. Then one can find a real-analytic and symplectic transformation $$\Phi:\ W_{v-2\hat v, s-2\hat s,\varepsilon-2\hat\varepsilon}\to W_{v, s,\varepsilon}$$ defined by the time-one flow\footnote{The time-one flow generated by $\phi$ is defined as the differential operator $$X_\phi^1:=\sum_{k=0}^\infty \frac{{\cal L}_\phi^k}{k!}$$ where ${\cal L}^0_\phi f:=f$ and ${\cal L}^k_\phi f:=\big\{\phi, {\cal L}^{k-1}_\phi f\big\}$, with $k=1,\ 2,\ \cdots$. } $X_\phi^1f:=f\circ\Phi$ of a suitable $\phi$ verifying $$\|\phi\|_{v,s,\varepsilon}\le \frac{\|f\|_{v,s,\varepsilon}}{\alpha_2}$$ such that $$H_+:=H\circ\Phi=h+g+\P_{{\mathbb L}} T_Kf+f_+$$ and, moreover, the following bounds hold \beqano \|f_+\|_{v-2\hat v, s-2\hat s, \varepsilon-2\hat\varepsilon}&\le& \big(1-\frac{c_m}{\alpha_2 \d}\|f\|_{v, s,\varepsilon}\big)^{-1}\Big[\frac{c_m}{\alpha_2 \d}\|f\|_{v, s,\varepsilon}^2 \nonumber\\
    &&+\max\left\{e^{-K\hat s/2}\,,\ \big(\frac{\varepsilon-\hat\varepsilon}{\varepsilon}\big)^{K/2}\right\}\|f\|_{v, s,\varepsilon}+ \|\big\{\phi, g \big\}\|_{v-\hat v, s-\hat s, \varepsilon-\hat\varepsilon}\Big] \eeqano Finally, for any real-analytic function $F$ on $W _{v, s, \varepsilon}$, \beqa{F-Id} \|F\circ\Phi-F\|_{v-2\hat v, s-2\hat s, \varepsilon-2\hat\varepsilon}&\le&\frac{\|\{\phi, F\}\|_{v-\hat v, s-\hat s, \varepsilon-\hat\varepsilon}}{\displaystyle  1-\frac{c_m\|f\|_{v, s,\varepsilon}}{\alpha_2 \d}}. \eeqa
\end{lemma}
{\bf Sketch of proof}  Lemma~\ref{itera1} is a straightforward generalization of~\cite[Iterative Lemma]{poschel93}. 
To obtain such generalization, just replace the norm defined in ~\cite[Section 1]{poschel93} with the norm~\equ{norm}, where
\beqa{exp f} f=\sum_{(k,\a,\b)\in \integer^{n}\times {\mathbb N}^\ell\times {\mathbb N}^\ell\atop{\a_i\ne \b_i\forall i}
  }f_{k, \a,\b}(I)e^{ik\cdot\varphi} \left(\frac{p-\ii q}{\sqrt 2}\right)^\a  \left(\frac{p+\ii q}{\ii\sqrt 2}\right)^\b\ ,\eeqa
and bound the ``ultraviolet remainders'', namely the norm of the functions whose expansion~\equ{exp f} includes only terms with 
$|(k, \alpha-\beta)|_1>K$, 
 as follows. Observe that, if $|(k, \alpha-\beta)|_1>K$, then either $|k|_1>K/2$  or $|\alpha-\beta|_1>K/2$. In the latter case, a fortiori, $|\alpha|_1+|\beta|_1\ge |\alpha-\beta|_1>K/2$. Then  we have, for such functions,
$\|f\|_{r, s-\hat s, \varepsilon-\hat\varepsilon}\le \max\left\{e^{-K\hat s/2}\,,\ \left(\frac{\varepsilon-\hat\varepsilon}{\varepsilon}\right)^{K/2}\right\}\|f\|_{r, s, \varepsilon}$. Other details are omitted.

\vskip.1in
\noindent
{\bf Proof of Proposition~\ref{averageELLIPTIC}} Let \beqano r_1:=r_0-2 \hat r_0,\qquad s_1:=s_0-2\hat s_0,\quad \varepsilon_1:=\varepsilon_0-2\hat\varepsilon_0\,. \eeqano By Lemma~\ref{itera1}, we find a canonical transformation $\Phi_1=X_{\phi_1}$ which is real-analytic on $W_{ r_1, s_1, \varepsilon_1}$ and conjugates $H=H_0$ to $H_1=H_0\circ\Phi_1=h+g_1+f_1$, where $g_1=\P_{{\mathbb L}}T_K f_0$ and

\beqano \|f_1\|_{v_1,s_1,\varepsilon_1}&\le&(1-\frac{c_mE_0}{\alpha_2 \d_0})^{-1}\Big[\frac{c_mE_0}{\alpha_2 \d_0}+\max\left\{e^{-K\hat s_0/2}\,,\ \big(\frac{\varepsilon_0-\hat\varepsilon_0}{\varepsilon_0}\big)^{K/2}\right\}\Big]E_0\nonumber\\
&\le&
2\Big[\frac{c_mE_0}{\alpha_2 \d_0}+e^{-K\hat \s_0/2}\Big]E_0
\eeqano 
having used 
$$\big(\frac{\varepsilon_0-\hat\varepsilon_0}{\varepsilon_0}\big)^{K/2}=e^{\frac{K}{2}\log \big(1-\frac{\hat\varepsilon_0}{\varepsilon_0}\big)}\le e^{-\frac{K}{2}\frac{\hat\varepsilon_0}{\varepsilon_0}}\,.$$
We now focus on the case
$$\frac{c_mE_0}{\alpha_2 \d_0}<e^{-K\hat \s_0/2}$$
otherwise the lemma is\footnote{Indeed, in such case,
$$ \|f_1\|_{v_1,s_1,\varepsilon_1}\le 4 e^{-K\hat \s_0/2}\le e^{-K\hat \s_0/4}$$
because $e^{-K\hat \s_0/4}\le \frac{1}{4}$ having chosen $K\hat \s_0\ge 8\log 2$.
} proved. Then we have
\beqano \|f_1\|_{v_1,s_1,\varepsilon_1}&\le&
4\frac{c_mE^2_0}{\alpha_2 \d_0}=:E_1\,.
\eeqano
Note that $$E_1<\frac{E_0}{4}\,.$$
Assume now that, for some $j\ge 1$, it is $H_j=H_{j-1}\circ \Phi_j=h+g_j+f_j$, where
\beqa{induction} g_j=\sum_{h=0}^{j-1}\P_{\mathbb L}T_K f_{h}\,,\qquad  \|f_j\|_{v_j,s_j,\varepsilon_j}&\le& E_j\le 
\min\left\{\frac{E_0}{4^j}\,,\ 4\frac{c_mE^2_0}{\alpha_2 \d_0}\right\}\,.
\eeqa
We have just proved this is true when $j=1$. Let $L:=\left[\frac{K\hat\sigma_0}{8\log 2}\right]$.
We prove that~\equ{induction} is true for $j+1$, for all $1\le j\le L$.
Let
\beqano\hat r_j:=\frac{\hat r_0}{L}\,,\quad \hat s_j:=\frac{\hat s_0}{L}\,,\quad \hat \varepsilon_j:=\frac{\hat \varepsilon_0}{L}\quad {\rm hence}\quad \delta_j=\frac{\delta_0}{L^2}\quad \forall\ 1\le j\le L\,.\eeqano
Note that, for all $1\le j\le L$, it is $\hat r_j<\frac{r_j}{2}$: $$r_j=r_1-2(j-1)\frac{\hat r_0}{L}\ge r_1-2(1-1/L)\hat r_0=r_0-4\hat r_0+2\hat r_j>2\hat r_j\,.$$
Similarly, $\hat s_j<\frac{s_j}{2}$, $\hat \varepsilon_j<\frac{\varepsilon_j}{2}$. 
Let then
$$ r_{j+1}=r_j-2\frac{\hat r_0}{L}\,,\quad s_{j+1}=s_j-2\frac{\hat s_0}{L}\,,\quad \varepsilon_{j+1}=\varepsilon_j-2\frac{\hat \varepsilon_0}{L}$$
so that $r_j=r_1-2(j-1)\frac{\hat r_0}{L}$, etc, for all $1\le j\le L$.
Then
\beqa{Ej small}c_m\frac{E_j}{\alpha_2 \delta_j}\le 4\frac{c^2_0E^2_0}{\alpha^2_2 \d^2_0}L^2<\frac{1}{16}\eeqa
and Lemma~\ref{itera1} applies again, and $H_j$ can be conjugated to $H_{j+1}=H_j\circ\Phi_{j+1}=h+g_{j+1}+f_{j+1}$, with
 \beqano 
 g_{j+1}&=&g_j+\P_{{\mathbb L}} T_Kf_j=\sum_{h=0}^{j}\P_{\mathbb L}T_K f_{h}\nonumber\\
 \|f_{j+1}\|_{r_{j+1}, s_{j+1}, \varepsilon_{j+1}}&\le& \big(1-\frac{c_m}{\alpha_2 \d_j}E_j\big)^{-1}\Big[\frac{c_m}{\alpha_2 \d_j}E_j^2 \nonumber\\
    &&+\max\left\{e^{-K\hat s_j/2}\,,\ \big(\frac{\varepsilon_j-\hat\varepsilon_j}{\varepsilon_j}\big)^{K/2}\right\}E_j+ \|\big\{\phi_j, g_j \big\}\|_{r_{j}-\hat r_j, s_j-\hat s_j, \varepsilon_j-\hat\varepsilon_j}\Big] \eeqano
    To bound the right hand side of the latter expression, we use~\equ{Ej small}
    and observe that    
\beqano
&& 
e^{-K\widehat s_j/2}=e^{-\frac{K}{2L}\widehat s_0}\le \frac{1}{16}\nonumber\\
&& \left(\frac{\varepsilon_j-\hat\varepsilon_j}{\varepsilon_j}\right)^{K/2}=\left(1-\frac{\frac{\widehat\varepsilon_0}{L}}{\varepsilon_1-2(j-1)\frac{\widehat\varepsilon_0}{L}}\right)^{K/2}\le \left(1-\frac{{\widehat\varepsilon_0}}{\varepsilon_1L}\right)^{K/2}\le e^{-\frac{K{\widehat\varepsilon_0}}{2\varepsilon_1L}}\le \frac{1}{16}\eeqano
having used 
$e^{-\frac{K\widehat s_0}{2L}}\le e^{-\frac{K\hat\sigma_0}{2L}}$, 
$e^{-\frac{K{\widehat\varepsilon_0}}{2\varepsilon_1L}}\le e^{-\frac{K{\widehat\varepsilon_0}}{2\varepsilon_0L}}\le e^{-\frac{K\hat\sigma_0}{2L}}$
 and $L\le\frac{K\hat\sigma_0}{8\log 2}$.  
Moreover, writing
$$g_j=\P_{\mathbb L}T_K f_{0}+{\mathbb 1}_{j\ge 2}\sum_{h=1}^{j-1}\P_{\mathbb L}T_K f_{h}=:f_{0}^{{\mathbb L}, K}+f_{j-1}^{{\mathbb L}, K}$$
with $f_{0}^{{\mathbb L}, K}$ real--analytic on $W_{r_0, s_0, \varepsilon_0}$, while $f_{j-1}^{{\mathbb L}, K}$ real--analytic on $W_{r_{j-1}, s_{j-1}, \varepsilon_{j-1}}$ and verifying
$$\|f_{0}^{{\mathbb L}, K}\|_{r_0, s_0, \varepsilon_0}\le E_0\,,\quad \|f_{j-1}^{{\mathbb L}, K}\|_{r_{j-1}, s_{j-1}, \varepsilon_{j-1}}\le \sum_{h=1}^{j-1}\frac{E_1}{4^{j-1}}\le \frac{4}{3}E_1$$
we get
\beqano
\|\big\{\phi_j, g_j \big\}\|_{r_{j}-\hat r_j, s_j-\hat s_j, \varepsilon_j-\hat\varepsilon_j}&\le& \frac{c_mL}{\alpha_2\delta_0}E_0{E_j}+\frac{4}{3}\frac{c_mL^2}{\alpha_2\delta_0}E_1{E_j}\nonumber\\
&\le&\left(\frac{c_mL}{\alpha_2\delta_0}E_0+\frac{16}{3}\frac{c^2_mL^2}{\alpha^2_2\delta^2_0}E_0^2\right){E_j}\nonumber\\
&\le&\left(\frac{1}{32}+\frac{1}{32}\right){E_j}=\frac{E_j}{16}
\eeqano
Collecting all such bounds we get
$$E_{j+1}\le \frac{16}{15}\frac{3}{16}E_j<\frac{E_j}{4}\,.$$
The inductive claim  $j\to j+1$ is thus proved, for all $1\le j\le L$. 
Letting now $\Phi_*:=\Phi_1\circ\cdots\circ \Phi_{L+1}$ and
$$H_*:=H_{L+1}=h+g_{L+1}+f_{L+1}=:h+g_{*}+f_{*}$$
$$ r_*:=r_{L+1}=r-4\hat r\,,\ s_*:=s_{L+1}=s-4\hat s\,,\ \varepsilon_*:=\varepsilon_{L+1}=\varepsilon-4\hat \varepsilon$$
and using $L+1> \frac{K\widehat \s_0}{8\log 2}$, we get
\beqano
&&\|f_*\|_{r_*, s_*, \varepsilon_*}\le \frac{E_0}{4^{L+1}}=e^{-2(L+1)\log 2}E_0< e^{-\frac{K\widehat \s_0}{4}}E_0\nonumber\\
&&\|g_*-\P_{\mathbb L}T_K f_{0}\|_{r_*, s_*, \varepsilon_*}\le \frac{4}{3}E_1<8\frac{c_mE^2_0}{\alpha_2 \d_0}
\eeqano
as claimed. The bounds~\equ{results2} are obtained from~\equ{F-Id}, by usual telescopic arguments. $\quad \square$

\paragraph*{Aknowledgments} G.P. is indebted to G. Gallavotti for many highlighting  discussions
.  Theorem~\ref{stable toriREF} is part of the MSc thesis of X.L. This research is funded by the PRIN project  "New frontiers of Celestial Mechanics: theory and applications" and has been developed under the auspices of INdAM and GNFM.\\ The authors confirm that the manuscript has no associated data.
\paragraph*{Author Contribution} GP conceptualized and developed the work and wrote the main manuscript. XL wrote a draft of Theorem~\ref{stable toriREF}, which is part of his MSc thesis at the University of Padova. All authors reviewed the manuscript.

\addcontentsline{toc}{section}{References}
 \bibliographystyle{plain}

\begin{thebibliography}{10}

\bibitem{arnold64}
V.~I. Arnold.
\newblock Instability of dynamical systems with many degrees of freedom.
\newblock {\em Dokl. Akad. Nauk SSSR}, 156:9--12, 1964.

\bibitem{arnold63}
V.I. Arnold.
\newblock {S}mall denominators and problems of stability of motion in classical
  and celestial mechanics.
\newblock {\em Russian Math. Surveys}, 18(6):85--191, 1963.

\bibitem{bonettoGGM1998}
F.~{B}onetto, G.~{G}allavotti, G.~{G}entile, and V.~{M}astropietro.
\newblock Lindstedt series, ultraviolet divergences and {Moser}'s theorem.
\newblock {\em Ann. Sc. Norm. Super. Pisa, Cl. Sci., IV. Ser.}, 26(3):545--593,
  1998.

\bibitem{cellettiC98}
A.~Celletti and L.~Chierchia.
\newblock Construction of stable periodic orbits for the spin-orbit problem of
  celestial mechanics.
\newblock {\em Regul. Chaotic Dyn.}, 3(3):107--121, 1998.
\newblock J. Moser at 70 (Russian).

\bibitem{chierchiaG94}
L.~Chierchia and G.~Gallavotti.
\newblock Drift and diffusion in phase space.
\newblock {\em Ann. Inst. H. Poincar\'e Phys. Th\'eor.}, 60(1):144, 1994.

\bibitem{chierchiaPi10}
L.~Chierchia and G.~Pinzari.
\newblock Properly--degenerate {KAM} theory (following {V}.{I}. {A}rnold).
\newblock {\em Discrete Contin. Dyn. Syst. Ser. S}, 3(4):545--578, 2010.

\bibitem{chierchiaPi11b}
L.~Chierchia and G.~Pinzari.
\newblock The planetary {$N$}-body problem: symplectic foliation, reductions
  and invariant tori.
\newblock {\em Invent. Math.}, 186(1):1--77, 2011.

\bibitem{chierchiaPr2019}
L.~{C}hierchia and M.~{P}rocesi.
\newblock {\em Kolmogorov-Arnold-Moser (KAM) Theory for Finite and Infinite
  Dimensional Systems}, pages 1--45.
\newblock Springer Berlin Heidelberg, Berlin, Heidelberg, 2019.

\bibitem{clarkeFG22}
A.~{C}larke{,} J.~F\'ejoz{,} and M.~Guardia.
\newblock Why are inner planets not inclided?
\newblock {\em ArXiv: 2210.11311v1}, 2022.

\bibitem{delshamsKDRS2019}
A.~{D}elshams, V.~{K}aloshin, A.~de~la {R}osa, and T.~M. {S}eara.
\newblock Global instability in the restricted planar elliptic three body
  problem.
\newblock {\em Commun. Math. Phys.}, 366(3):1173--1228, 2019.

\bibitem{deprit83}
A.~Deprit.
\newblock Elimination of the nodes in problems of {$n$} bodies.
\newblock {\em Celestial Mech.}, 30(2):181--195, 1983.

\bibitem{fejoz04}
J.~F{\'e}joz.
\newblock D{\'e}monstration du `th{\'e}or{\`e}me d'{A}rnold' sur la
  stabilit{\'e} du syst{\`e}me plan{\'e}taire (d'apr{\`e}s {H}erman).
\newblock {\em Ergodic Theory Dynam. Systems}, 24(5):1521--1582, 2004.

\bibitem{FGKR14}
J.~F{\'e}joz, M.~Guardia, V.~Kaloshin, and P.~Roldan.
\newblock {K}irkwood gaps and diffusion along mean motion resonances in the
  restricted planar three body problem.
\newblock {\em J. Eur. Math. Soc.}, 2014.

\bibitem{gallavottiG1995}
G.~{G}entile G.~{G}allavotti.
\newblock Majorant convergence for twistless {KAM} tori.
\newblock {\em Ergodic Theory Dyn. Syst.}, 15(5):587--869, 1995.

\bibitem{gallavotti86}
G.~Gallavotti.
\newblock Quasi-integrable mechanical systems.
\newblock In {\em Ph\'enom\`enes critiques, syst\`emes al\'eatoires, th\'eories
  de jauge, {P}art {I}, {II} ({L}es {H}ouches, 1984)}, pages 539--624.
  North-Holland, Amsterdam, 1986.

\bibitem{gallavotti1994}
G.~{G}allavotti.
\newblock Twistless {KAM} tori.
\newblock {\em Commun. Math. Phys.}, 164(1):145--156, 1994.

\bibitem{guzzoEP2020}
M.~Guzzo, C.~Efthymiopoulos, and R.~I. Paez.
\newblock Semi-analytic computations of the speed of {Arnold} diffusion along
  single resonances in a priori stable {Hamiltonian} systems.
\newblock {\em J. Nonlinear Sci.}, 30(3):851--901, 2020.

\bibitem{jacobi1842}
C.~G.~J. Jacobi.
\newblock Sur l'{\'e}limination des noeuds dans le probl{\`e}me des trois
  corps.
\newblock {\em Astronomische Nachrichten}, Bd XX:81--102, 1842.

\bibitem{laskarR95}
J.~Laskar and P.~Robutel.
\newblock Stability of the planetary three-body problem. {I}. {E}xpansion of
  the planetary {H}amiltonian.
\newblock {\em Celestial Mech. Dynam. Astronom.}, 62(3):193--217, 1995.

\bibitem{pinzari18}
G.~{Pinzari}.
\newblock {On the co-existence of maximal and whiskered tori in the planetary
  three-body problem}.
\newblock {\em Journal of Mathematical Physics}, 59(5):052701, May 2018.

\bibitem{pinzari18a}
G.~Pinzari.
\newblock Perihelia reduction and global {K}olmogorov tori in the planetary
  problem.
\newblock {\em Mem. Amer. Math. Soc.}, 255(1218), 2018.

\bibitem{poschel93}
J.~P{\"o}schel.
\newblock Nekhoroshev estimates for quasi-convex {H}amiltonian systems.
\newblock {\em Math. Z.}, 213(2):187--216, 1993.

\bibitem{radau1868}
R.~Radau.
\newblock Sur une transformation des {\'e}quations diff{\'e}rentielles de la
  dynamique.
\newblock {\em Ann. Sci. Ec. Norm. Sup.}, 5:311--375, 1868.

\bibitem{Russmann:2001}
H.~R{\"u}ssmann.
\newblock Invariant tori in non-degenerate nearly integrable {H}amiltonian
  systems.
\newblock {\em Regul. Chaotic Dyn.}, 6(2):119--204, 2001.

\end{thebibliography}
\def\cprime{$'$} \def\cprime{$'$}

\end{document}